\documentclass[letterpaper,11pt]{article}

\usepackage{amsmath}
\usepackage{amsthm}

\usepackage {amssymb}
\usepackage {stmaryrd}
\usepackage[final]{graphicx}
\usepackage[all]{xy}

\newtheorem{theorem}{Theorem}
\newtheorem{cor}[theorem]{Corollary}
\newtheorem{lemma}[theorem]{Lemma}
\newtheorem{definition}[theorem]{Definition}

\newtheorem{remark}{Remark}

\title{\bf Gr\"obner bases, monomial group actions, and the Cox rings of Del Pezzo surfaces}
\author{Mike Stillman, Damiano Testa, Mauricio Velasco}
\date{\begin{small}
Department of Mathematics, Cornell University, Ithaca, NY 14853, USA\end{small}}
\begin{document}
\maketitle
\abstract{
We introduce the notion of monomial group action and study some
of its consequences for Gr\"obner basis theory. As an application we
prove a conjecture of V. Batyrev and O. Popov describing the Cox rings
of Del Pezzo surfaces (of degree $\geq 3$) as quotients of a
polynomial ring by an ideal generated by quadrics.}
\section{Introduction}
The notion of homogeneous coordinate ring was introduced by David Cox in~\cite{COX} aiming to generalize to arbitrary toric varieties the relationship between $\mathbb{P}^n$ and $k[x_0,\dots,x_n]$. Cox's construction assigns to every toric variety $T$, a multigraded {\it polynomial ring} $R$ (and an ideal) such that:
\begin{enumerate}
\item{$T$ can be recovered as a suitable quotient of ${\rm Spec} (R)$ by the action of a torus;}
\item{Modules over $R$ correspond to sheaves on $T$.}
\end{enumerate}
This construction was generalized by Keel and Hu in~\cite{HK} where the authors introduce {\it Cox rings}, the homogeneous coordinate rings of a much larger class of varieties. The authors show that finite generation of this ring is of fundamental importance for the birational geometry of the variety (in particular, it ensures that the Mori program can be carried out for any divisor see~\cite{HK} prop. 1.11).\\
Moreover Keel and Hu prove that toric varieties are the only algebraic varieties whose Cox rings are polynomial rings, thus raising the question of which kinds of finitely generated $k$-algebras arise as Cox rings of non-toric varieties.\\
Probably the most important such example are the Cox rings of Del Pezzo surfaces of degree at most five. These rings were studied for the first time by Batyrev and Popov in~\cite{BP}, where the authors show that they are Gorenstein $k$-algebras whose generators are in bijection with the $(-1)$-curves on the surfaces. Moreover, they conjecture that these rings are quadratic algebras.\\
\\
This paper is a case study of the Cox rings of Del Pezzo surfaces, specifically of the ideals $C$ which define them as quotients of the polynomial rings $k[E]$ (with one variable for each exceptional curve).\\ 
\\
The groups of symmetries of the configuration of exceptional curves play a fundamental role in our study. We show that, although the action of this Weyl group on $k[E]$ does not fix the ideal $C$, it can be rediscovered as symmetries of the Gr\"obner fan of $C$.\\
This weaker form of symmetry is sufficient to characterize the monomial initial ideals of $C$ in terms of very few values of their multigraded Hilbert Series. As an application of the techniques developed we prove Batyrev and Popov's conjecture (for surfaces of degree at least three) providing explicit generators for $C$.\\
\\
The material is organized as follows: 
\begin{itemize}
\item{Section 2 contains background material on Del Pezzo surfaces and their Weyl groups.} 
\item{Section 3 contains the definition Cox rings and the results of Batyrev and Popov used throughout the rest of the paper.}
\item{In Section 4 we describe the degree 2 part of the ideals $C$ which define the Cox rings as quotients of polynomial rings.}
\item{In Section 5 we introduce the notion of {\it monomial group action}, and study its consequences for Gr\"obner basis theory. In particular we show that, if $G$ acts monomially on an  ideal $I$, then it acts by symmetries on its Gr\"obner fan (and in particular on its tropical variety). Moreover we show that, for a general Del Pezzo surface, the corresponding Weyl group acts monomially on the ideal defining its Cox ring.}
\item{In Section 6 we study the problem of characterizing the monomial initial ideals of a homogeneous ideal in the presence of a monomial group action compatible with the grading. As an application we characterize the monomial initial ideals of $C$ in terms of very few values of their Picard-graded Hilbert Series.}
\item{Section 7 contains the proof of Batyrev and Popov's conjecture (for surfaces of degree at least three). Moreover, we show that for degree al least 4 the Cox rings of Del Pezzo surfaces are Koszul algebras.}
\item{In Section 8 we reduce the problem of finding quadratic Gr\"obner bases for $C$ to a combinatorial problem about the edge ideals of the graphs of exceptional curves.}
\item{Section 9 is an Appendix containing tables and calculations used in Section 7.}  
\end{itemize}

\section{Del Pezzo surfaces}
This section contains the required background on Del Pezzo surfaces and their Weyl groups, and introduces terminology that will be used throughout the rest of the paper.
\begin{definition} A collection of $r\leq 8$ points in $\mathbb{P}^2$ is said to be in general position if no three are on a line, no six are on a conic and any cubic containing eight points is smooth at each of them.\end{definition}
\begin{definition} A Del Pezzo surface $X_r$ is the blowup of $\mathbb{P}^2$ at $r\leq 8$ general points. The degree of $X_r$ is $9-r$.\end{definition}

\begin{remark}
Normally the definition of Del Pezzo surfaces includes $\mathbb{P}^1 \times \mathbb{P}^1$ of degree 8, but we concentrate on Del Pezzo surfaces of degree at most 5.
\end{remark}

Since the automorphism group of $\mathbb{P}^2$ carries any four general points to the four standard ones, there is essentially one Del Pezzo surface $X_r$ for $r\leq 4$. In constrast, there are infinitely many nonisomorphic Del Pezzo surfaces $X_r$ for each $r\geq 5$.\\
\\
From the description of Del Pezzo surfaces as blow ups of $\mathbb{P}^2$ it follows immediately that the Picard group of $X_r$ is isomorphic to $\mathbb{Z}^{r+1}$. A natural basis is given by:
\begin{itemize}
\item{The pullback of the class of a line in $\mathbb{P}^2$, denoted by  $\ell$;}
\item{The exceptional divisors of the blow up $e_1,\dots, e_r$.}
\end{itemize}
In terms of this basis, the intersection form is given by 
\[e_i\cdot e_j=-\delta_{ij} \text{ , } \ell\cdot \ell=1\text{ and } \ell\cdot e_j=0\]
Moreover the canonical divisor in $X_r$ is $K=-3\ell+e_1+$\dots$+e_r$\\
\\
For $r\leq 6$, the linear system associated to the anticanonical divisor determines an embedding of $X_r$ as a surface of degree $9-r$ in $\mathbb{P}^{9-r}$.\\ 
The best known examples of Del Pezzo surfaces are $\mathbb{P}^2$, the $X_5$, embedded by $|-K|$ as complete intersections of two quadrics in $\mathbb{P}^4$, and the $X_6$, which correspond to smooth cubic surfaces in $\mathbb{P}^3$.\\
\\
Del Pezzo surfaces contain a very special collection of rational curves (for $r=6$ these curves are the 27 lines on the cubic).
\begin{definition}
An exceptional curve (also $(-1)$-curve)  $C$ is a curve whose class in ${\rm Pic}(X_r)$ satisfies:
\[ K\cdot C = -1 \text{ ~ and ~ } C^2=-1\]
\end{definition}
Each Del Pezzo surface contains finitely many exceptional curves which are classified (for $r\leq 6$) by the following table (see~\cite{Man} for details)
\begin{center}
\begin{tabular}{c|ccc|c}
Number of blown up points & $4$ & $5$ & $6$ & \text{Class in ${\rm Pic}(X_r)$}\\
\hline
\hline
Exceptional divisors $e_i$ & $4$ & $5$ & $6$ & $e_i$\\
Lines through pairs of points $f_{ij}$ & $6$ & $10$ & $15$ & $\ell-e_i-e_j$\\
Conics through five points $g_i$ & $0$ & $1$ & $6$ & $2\ell-\sum_{k\neq i}e_k$\\
\hline
\hline
Number of exceptional curves & $10$ & $16$ & $27$\\
\end{tabular}
\end{center}
\bigskip
The configuration of the exceptional curves on the surface is better visualized by a graph (with multiple edges for $r\geq 7$).\\
\begin{definition} The graph of $(-1)$-curves is the graph with one vertex for each exceptional curve and  $E_i\cdot E_j$ edges between edges $E_i$ and $E_j$ for $i\neq j$.  We denote this graph by $L_r$.\end{definition} 
The graphs $L_4$ and $L_5$ (the Petersen and Clebsch graphs respectively) are shown in the figure. Note that the configuration of lines is independent of the coordinates of the blown up points.
\begin{center}
\scalebox {1.1}{
\includegraphics{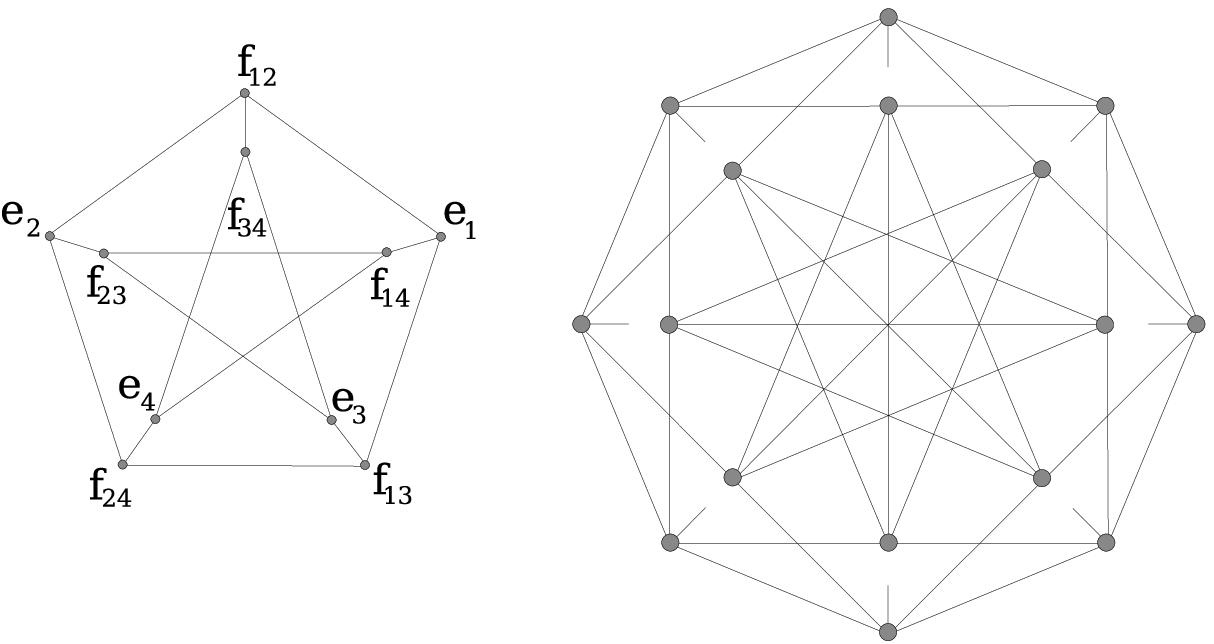}
}
\end{center}
\subsection{Symmetries}
For each $r\geq 2$ there is a Weyl group $W_r$ which acts on ${\rm Pic}(X_r)$ by automorphisms which preserve the intersection form.
\begin{center}
\begin{tabular}{c|c|c}
N. of Blown up points & Root system & size of $W_r$\\
\hline
\hline
$4$ & $A_4$ & $120$\\
$5$ & $D_5$ & $1920$\\
$6$ & $E_6$ & $51840$\\
\hline
\end{tabular}
\end{center}
\bigskip
More concretely, $W_r$ is the subgroup of ${\rm Aut}({\rm Pic}(X_r))$ generated by the permutations of the classes of the exceptional divisors $e_i$ and (for $r\geq 3$) by the additional Cremona element $\sigma$ given by $\sigma(\ell)=2\ell-e_1-e_2-e_3$, $\sigma(e_1)=\ell-e_2-e_3$, $\sigma(e_2)=\ell-e_1-e_3$, $\sigma(e_3)=\ell-e_1-e_2$ and $\sigma(e_i)=e_i$ for $i\not\in \{1,2,3\}$.\\
\\
The elements of $W_r$ preserve the intersection form and fix the canonical divisor $K$. As a result they permute the classes of $(-1)$-curves (since they fix the equations that define them in ${\rm Pic}(X_r)$) and induce automorphisms of the graphs of $(-1)$-curves (since these permutations preserve intersection numbers). The transitivity of this action on vertices and edges explains the striking symmetry of the graphs.
\section{Cox rings}
The following definition was proposed by Hu and Keel in~\cite{HK}
\begin{definition} Let $X$ be a projective variety with $N^1(X)={\rm Pic}(X)_{\mathbb{Q}}$ and let $L_1,\dots, L_k$ be line bundles which are a basis for the torsion free part of the Picard group and whose affine hull contains $\overline{NE}^1(X)$. A Cox ring for $X$ is the ring 
\[CR(X,L_1,\dots,L_k)=\bigoplus_{(m_1,\dots,m_k)\in\mathbb{Z}^k}H^0\bigl( L_1^{\otimes m_1}\otimes\dots\otimes L_k^{\otimes m_k}\bigr)\]
\end{definition}
\noindent Note that the isomorphism type of this ring is independent of the choice of basis (see~\cite{HK} for details).\\
\\
For a Del Pezzo surface $X_r$, we choose the following basis of ${\rm Pic}(X_r)$: 
\begin{itemize}
\item{The $r$ line bundles ${\cal O}[e_i]$ corresponding to the exceptional divisors $e_i$;}
\item{ ${\cal O}[\ell]$ where $\ell$ is the pullback of the line $z=0$ in $\mathbb{P}^2$.}
\end{itemize}
\begin{definition} We denote by $Cox(X_r)$ the ring
\begin{eqnarray*}
Cox(X_r) & = & CR(X_r;\ell,e_1,\dots,e_r)= \\
& = & \hspace{-6pt} \bigoplus_{(m_0,\dots,m_r)\in\mathbb{Z}^r} \hspace{-6pt} 
H^0({\cal O}[m_0 \ell+m_1 e_1+\dots+m_r e_r])
\end{eqnarray*}
\end{definition}
\noindent It is obvious from the definition that $Cox(X_r)$ is a ${\rm Pic}(X_r)$-graded integral domain. Moreover this ring admits a coarser $\mathbb{Z}$-grading given by
\[Cox(X_r)_n=\hspace{-10pt} \bigoplus_{\{D\in {\rm Pic}(X_r):-K\cdot D=n\}} \hspace{-15pt} Cox(X_r)_D \hspace{15pt} \text{ for $n\in\mathbb{Z}$}\]
Note that the above grading is nonnegative.\\
\\
The Cox rings of Del Pezzo surfaces were studied for the first time by Batyrev and Popov in~\cite{BP}  where they show the following fundamental result:
\begin{theorem} \label{Thm:BP} For $3\leq r \leq 7$ the ring $Cox(X_r)$ is generated by the global sections of invertible sheaves defining the exceptional curves.\end{theorem}
\noindent In particular, effective divisor classes can be written as sums of classes of exceptional curves.
\section{Cox rings of Del Pezzo surfaces as quotients of polynomial rings}
\noindent For $r<4$ the Del Pezzo surface $X_r$ is a toric variety and its Cox ring is a polynomial ring (see~\cite{COX}). In this section we set up the notation necessary to describe $Cox(X_r)$ for $4\leq r\leq 7$.
\begin{definition} Let $E_r$ be the set of Picard classes of exceptional curves in $X_r$. Let $k[E_r]$ be the ${\rm Pic}(X_r)$-graded polynomial ring  obtained by letting $\deg([c])=[c]$. \end{definition}
\noindent
For clarity we use the symbols $e_i, f_{ij}, g_i$ (in correspondence with the exceptional divisors, the strict transforms of lines through pairs of points and the strict transforms of conics through five points resp.) as variables in $k[E_r]$.\\
\\
In this notation, Theorem~\ref{Thm:BP} shows that, for every choice of nonzero global sections $s_D\in H^0 ( \mathcal {O} [D])$ with $D=m_0\ell+m_1e_1+\dots+m_re_r$ such that $[D]$ is the class of an exceptional curve, the map $\phi:k[E_r]\rightarrow Cox(X_r)$ which sends each variable to the corresponding section is a ${\rm Pic}(X_r)$-graded surjective homomorphism.\\ 
Note that $s_D$ is determined by $D$ only up to multiplication by a nonzero constant so there are many possible maps $\phi$. 
\begin{definition}\label{Def:C_r}
We denote by $C_r$ the ${\rm Pic}(X_r)-$homogeneous prime ideal $ker(\phi)$ for some choice of
sections $s_D$.  In particular $Cox(X_r)\cong k[E_r]/C_r$.
\end{definition}
\noindent Note that $C_r$ depends on the choice of sections. However, multipliying the variables by constants is an automorphism of $k[E_r]$ which carries any choice of $C_r$ to any other.\\
\\
The images of monomials of $k[E_r]$ in $Cox(X_r)$ are of particular importance, we call their multiples by a nonzero constant {\it distinguished global sections}.
\begin{definition} A section of a bundle $D$ is distinguished if it is supported in a union of exceptional curves on $X_r$.\end{definition}
\noindent Note that the linear dependencies between distinguished global sections generate the ideal $C_r$.  Describing the ideal $C_r$ explicitly is one of the objectives of this paper, we begin by describing its (coarse) degree $2$ part (as in~\cite{BP}).
\begin{definition} A divisor class $C$ on $X_r$ is called a conic if it satisfies
\[-K\cdot C=2\text{ ~ and ~ } C^2=0\]\end{definition}
\noindent It is an easy consequence of Riemann-Roch and the adjunction formula that if $C$ is a conic then the linear system $|C|$ is base-point free and induces a morphism $X_r\rightarrow \mathbb{P}^1$ which is a conic bundle.\\
Moreover every such divisor has exactly $r-1$ distinguished global sections and any set of three of them are linearly dependent (3 vectors in a 2 dimensional vector space). As a result every conic provides $r-3$ linearly independent elements of $C_r$. 
\begin{definition}
We denote by $Q_r\subset C_r$ the ideal generated by the linear dependencies among distinguished global sections of conics $D$.
\end{definition}
\noindent Note that $-K\cdot D=2$ implies that  these relations are quadrics in $k[E_r]$ and is easy to see that they generate the degree 2 part of $C_r$ (since every $D$ with $-K\cdot D=2$ is either a conic or contains exactly one distinguished global section). \\
\\
The conic bundles $C$ on $X_r$ for $r\leq 6$ are described (up to permutation of the $e_i$'s) in the following table:\\
\begin{center}
 \begin{tabular}{c|ccc|c}
$r$ & $4$ & $5$ & $6$ & \text{Class in ${\rm Pic}(X_r)$}\\
\hline
\hline
& $4$ & $5$ & $6$ & $\ell-e_i$\\
& $1$ & $5$ & $15$ & $2\ell-e_1-e_2-e_3-e_4$\\
& $0$ & $0$ & $6$ & $3\ell-2e_1-e_2-e_3-e_4-e_5-e_6$\\
\hline
\hline
& $5$ & $10$ & $27$ & \text{Total number of conics}\\
& $5$ & $20$ & $81$ & \text{Total number of generators of $Q_r$}\\
\end{tabular}
\end{center}
\noindent The ideals $Q_r(p_1,\dots, p_r)$ can be generated in Macaulay2 using our package $CRDelPezzo$, which will be a part of the Macaulay2 distribution.
\\
{\bf Example.} For $r=4$, with blown up points $[1:0:0]$, $[0:1:0]$, $[0:0:1]$ and $[1:1:1]$ we have:
\begin{itemize}
\item{$k[E_4]=k[f_{12},\dots,f_{34},e_1,\dots ,e_4]$ graded by
$$ \deg(f_{ij})=\ell-\underline{e_i}-\underline{e_j} \text{ ,  } \deg(e_i)=\underline{e_i} $$}
\item{$Q_4$ is the ideal generated by
\[\begin{array}{c|c}
$Conic$ & $Generator$\\
\hline
2\ell-e_1-e_2-e_3-e_4 & f_{14}f_{23}-f_{12}f_{34} -f_{13}f_{24},\\
\ell-e_1 & e_2f_{12} -e_3f_{13}-e_4f_{14},\\
\ell-e_2 & e_1f_{12} - e_3f_{23}-e_4f_{24},\\
\ell-e_3 & e_1f_{13} -e_2f_{23}+e_4f_{34},\\
\ell-e_4 & e_1f_{14} -e_2f_{24}-e_3f_{34}\\
\end{array}\]
}
\end{itemize}
{\bf A conjecture of Batyrev and Popov.} In~\cite{BP} the authors conjecture that for every $4\leq r\leq 8$ and for every Del Pezzo surface $X_r$ the ideal $C_r$ is generated by quadrics  (i.e. $Q_r=C_r$).\\ 
\\
Batyrev and Popov observe that the equality $Q_r=C_r$ holds up to radical. Moreover they prove that $Q_4=C_4$ by observing that $k[E_4]/Q_4$ is the homogeneous coordinate ring of the Grassmannian $Gr(2,5)$ and hence an integral domain.\\
\\
We prove Batyrev and Popov's conjecture for $r=4,5$ and cubic surfaces without Eckart points in Section 7. Our proof does not depend on the equality between the radicals of these ideals. It is our hope that the methods developed here could be used to characterize the ideals of relations of other Cox rings which are known to be finitely generated $k$-algebras.\\

\section{Monomial group actions}

Throughout the rest of the section $R=k[x_1,\dots,x_n]$ denotes the ring of polynomials over a field $k$ with the standard grading and $G$ is a group acting on $R$ by permuting the variables; $I \subset R$ is a homogeneous ideal. For the necessary background on Gr\"obner basis see~\cite{GB}.

\begin{definition} For $h\in R$, ${\rm mon} (h)$ is the set of monomials of $R$ which appear with nonzero coefficient in $h$.
\end{definition}
\begin{definition}
The group $G$ acts monomially on $I$ up to degree $d$ if for every $h\in I$ of degree $\leq d$ 
and every $g\in G$ there is an element $h' \in I$ such that ${\rm mon} (h') = {\rm mon} \bigl( g (h) \bigr)$.
\end{definition}
Note that, if $G$ acts monomially on $I$ up to degree $d$ then it acts on the set $\{ {\rm mon}(h):h\in I\text{ and }\deg(h)\leq d\}$.
\begin{definition} Given a monomial order $\preceq$ and an element $g\in G$, let $\preceq_{g}$ be the monomial order given by
\[a\preceq_{g} b \Leftrightarrow g(a)\preceq g(b)\]  
\end{definition}
\begin{lemma}\label{lem:weightsmonomial} If $G$ acts monomially on $I$ up to degree $d$ and $in_{\preceq}(I)$ is generated in degree $\leq d$ then  $g(in_{\preceq}(I))=in_{\preceq_{g^{-1}}}(I)$ for any $g\in G$.\end{lemma}
\begin{proof} Let $S=\{s_1,\dots, s_k\}$ be a reduced $\preceq$ Gr\"obner basis for $I$ and let $g\in G$. By the monomiality of the action there is a set $S'=\{h_1,\dots, h_k\}\subset I$ such that ${\rm mon}(h_i)=g({\rm mon}(s_i))$. We show that $S'$ is a $\preceq_{g^{-1}}$ Gr\"obner basis. By definition of $\preceq_{g^{-1}}$, $in_{\preceq_{g^{-1}}}(h_i)=g(in_{\preceq}(p_i))$ so
\[in_{\preceq_{g^{-1}}}(I)\supseteq \bigl( in_{\preceq_{g^{-1}}}(h_i)\bigr)\supseteq g\bigl((in_{\preceq}(p_i))\bigr)=g(in_{\preceq}(I))\]
and all these ideals coincide since the first and the last have the same Hilbert function as $I$ ($G$ acts on $R$ by automorphisms of standard degree $0$). As a result $S'$ is a $\preceq_{g^{-1}}$ reduced Gr\"obner basis.
\end{proof}

In particular $G$ acts on the set of monomial initial ideals of $I$. We show that in fact this action extends to the Gr\"obner Fan of $I$.
\begin{definition} Given a weight vector $w=(w_1,\dots,w_n)$ and $g\in G$, let $g(w)=(w_{g(1)},\dots, w_{g(n)})$.\end{definition}
\begin{definition} For a weight vector $w$ let 
$$ C[w]=\{w'\in\mathbb{R}^n:\text{ } in_{w'}(I)=in_{w}(I)\}$$
and let $\preceq_{w}$ be the monomial order defined by 
\[a\preceq_{w} b \Leftrightarrow w\cdot a<w\cdot b \text{ or } w\cdot a=w\cdot b\text{ and }
a\preceq b\] 
where $\preceq$ is a fixed monomial term order and we have identified monomials with their exponent vectors.
\end{definition}
\noindent We denote by $\preceq_{gw}$ the monomial order $(\preceq_{w})_g$. Note that $\preceq_{gw}$ refines the preorder given by the weight $g(w)$.   
\begin{lemma} \label{lem:gfan} Assume that the largest degree of a generator in any monomial initial ideal of $I$ is $d$. If $G$ acts monomially on $I$ up to degree $d$ then the action of $G$ on weight vectors induces automorphisms of the Gr\"obner fan of the ideal $I$.\end{lemma}
\begin{proof}
Let $S$ be a reduced $\preceq_{w}$ Gr\"obner basis and let $S'$ be a reduced $\preceq_{g^{-1}w}$ Gr\"obner basis constructed as in Lemma~\ref{lem:weightsmonomial}.\\ 
By Proposition 2.3 in~\cite{GB}, 
\[C[g^{-1}(w)]=\{\eta\in \mathbb{R}^n: in_{\eta}(s')=in_{g^{-1}w}(s')\text{ for all }s'\in S'\}\]
By definition of $g^{-1}(w)$ and $S'$, $\eta \in C[g^{-1}(w)]$ if and only if $g(\eta)$ chooses the same leading forms of the generators of $S$ as the weight $w$ does.\end{proof} 
Similarly, under the above conditions $in_{w}(I)$ contains a monomial $m$ if and only if $in_{g(w)}(I)$ contains the monomial $g^{-1}(m)$. As a result
\begin{cor} Assume that the largest degree of a generator in any monomial initial ideal of $I$ is $d$. If $G$ acts monomially on $I$ up to degree $d$ then $G$ acts by automorphisms on the tropical variety ${\cal T}(I)$.\end{cor}
Now we show that the Weyl group $W_r$ acts monomially on $C_r$ for a general choice of blown up points. It should be remarked that the action of $W_r$ on $k[E_r]$ by permutation of the coordinates does not, in general, fix the ideals $C_r$ (even if we allow the permutations to be precomposed with diagonal matrices) so that monomiality is a way to recover the symmetries present  in the configuration of exceptional curves. These reappear as symmetries of the Gr\"obner fan of $C_r$.\\
\\
For clarity we denote by $C_r(p_1,\dots, p_r)$ the ideal of relations of the Cox ring of the Del Pezzo surface obtained by blowing up $\mathbb{P}^2$ at the points $p_1,\dots, p_r$.
\begin{definition} The group $W_r$ acts on $(\mathbb {P}^2)^r$ by birational automorphisms. We let the symmetric group on the $r$ indices act by permuting the coordinates and (if $r \geq 3$) we let the generator $T_{123}$ act by 
\begin{eqnarray*}
(\mathbb {P}^2)^r & \dashrightarrow & (\mathbb {P}^2)^r \\
(p_1 , \ldots , p_r) & \longmapsto & \bigl( q_1 , q_2 , q_3 , C(p_4) , \ldots , C(p_r) \bigr)
\end{eqnarray*}
where $C :  \mathbb {P}^2 \dashrightarrow \mathbb {P}^2$ is the Cremona transformation based at 
$p_1 , p_2 , p_3$ and $q_i$ is the image under $C$ of the line between $p_j$ and $p_k$, 
$\{ i , j , k \} = \{ 1 , 2 , 3 \}$.  \end{definition}
\noindent
Note that the action of $W_r$ restricts to automorphisms of the open subset of $(\mathbb{P}^2)^r$ consisting of points $(p_1 , \ldots , p_r)$ such that the blow up of $\mathbb {P}^2$ at $p_1 , \ldots , p_r$ is a del Pezzo surface. We denote this open subset by $U$.\\
\\
Let $\phi_{\lambda}$ be the automorphism of $k[E_r]$ obtained by multiplying each variable by a component of a vector $\lambda$ of nonzero constants.
\begin{lemma} \label{lem:WDelPezzos} For all $g\in W_r$ there exists a $\lambda$ such that 
$$ \phi_{\lambda}(C_r(g(p_1,\dots, p_r)) = g(C_r(p_1,\dots, p_r))$$
In particular $m_1,\dots, m_k$ are the monomials of an element of the ideal $C_r(p_1,\dots, p_r)$
if and only if 
$g(m_1),\dots, g(m_k)$ are those in an element of $C_r(g(p_1,\dots,p_r))$.\end{lemma}
\begin{proof} Recall from Definition~\ref{Def:C_r} that $C_r(p_1,\dots, p_r)$ depends on the choice of sections $s_D$ and that multiplying the variables by nonzero constants carries any choice of $C_r(p_1,\dots, p_r)$ to any other. Now, linear relations between distinguished global sections are intrinsic to the surface and do not depend on its particular presentation as a blow up of $\mathbb{P}^2$. 
Since the blow up of $\mathbb{P}^2$ at $(p_1,\dots, p_r)$ and at $g(p_1,\dots, p_r)$ are two presentations of the same Del pezzo surface the result follows.\end{proof}
\begin{theorem}\label{thm:monomialita}
For any integer $d$ there exists an open dense set of choices for $p_1,\dots, p_r$ 
such that the group $W_r$ acts monomially on the ideal $C_r(p_1,\dots, p_r)$ up to 
degree $d$. 
\end{theorem}
\begin{proof}
Denote by $M_{\leq d}$ the set of all subsets of monomials of $k[E_r]$ of coarse degree at most $d$.\\ For any $\underline m \in M_{\leq d}$ the locus $D_{\underline m}$ in $(\mathbb {P}^2)^r$ of points 
$(p_1 , \ldots , p_r)$ such that the monomials in $\underline m$ are linearly dependent in $k[E_r]/ C_r(p_1,\dots,p_r)$ is closed: it is the locus where the rank of ${\rm span} (\underline m) \rightarrow R / C_r$ is smaller than $|\underline m|$.  We partition the set $M_{\leq d}$ into two disjoint subsets $\mathcal G$ and $\mathcal S$, setting, for each 
$\underline m \in M_{\leq d}$, 
\begin{itemize}
\item $\underline m \in \mathcal G$ if $D_m$ is a proper closed subset of 
$(\mathbb {P}^2)^r$, 
\item $\underline m \in \mathcal S$ if $D_m = (\mathbb {P}^2)^r$.
\end{itemize}

The set $\mathcal G$ (and thus $\mathcal S$) is a union of orbits of $G$, since by Lemma~\ref{lem:WDelPezzos},  $D_{\underline m}$ does not contain $(p_1 , \ldots , p_r)$ if and only if $g \cdot D_{\underline m}$ does not contain $g \cdot (p_1 , \ldots , p_r)$.\\
Define 
$$ \mathcal B = \bigl( \mathbb {P}^2 \bigr)^r \setminus \mathop {\cup } \limits_{\underline m \in \mathcal S} D_{\underline m} $$
and observe that $\mathcal B$ is disjoint from $D_{\underline m}$, for all $\underline m \in \mathcal S$.\\
Let $(p_1 , \ldots , p_r) \in \mathcal B$ and let $X = Bl_{p_1 , \ldots , p_r} (\mathbb {P}^2)$.  We want to prove that 
the action of $G$ is monomial on $C_r(p_1,\dots,p_r)$.  We proceed by induction on the set ${\rm mon} (C_r) \cap M_{\leq d}$ (partially) 
ordered by inclusion.\\
Suppose that $f \in C_r$, $f \neq 0$ and ${\rm mon} (f)$ is minimal under inclusion in ${\rm mon}(C_r) \cap M_{\leq d}$.  
In particular $D_{{\rm mon} (f)} \cap \mathcal B \neq \emptyset$ (as it contains $(p_1 , \ldots , p_r)$) and therefore 
${\rm mon} (f) \in \mathcal G$.  Since the set $\mathcal G$ is $W_r-$invariant, for any $g \in G$ we have that 
$g \cdot {\rm mon} (f) \in \mathcal G$.  Thus $D_{g \cdot {\rm mon} (f)} \supset \mathcal B$, and therefore there is a non-zero 
relation $f^g \in C_r$ such that ${\rm mon} (f^g) \subset g \cdot {\rm mon} (f)$.  Note that typically $g \cdot f \not \in C_r$.  The 
same argument applied to $f^g$ and the group element $g^{-1}$ implies that there is a non-zero relation 
$(f^g)^{g^{-1}} \in C_r$ such that ${\rm mon} \bigl( (f^g)^{g^{-1}} \bigr) \subset g^{-1} \cdot g \cdot {\rm mon} (f) = {\rm mon} (f)$.  
By the minimality of ${\rm mon} (f)$ it follows that all the inclusions of monomials were in fact equalities and we conclude the proof of the base case of the induction.

Suppose that $f \in C_r$, $f \neq 0$ and that there is a non-zero $h \in C_r$ such that 
${\rm mon} (h) \subsetneq {\rm mon} (f)$.  Then we can find a constant $t$ such that $h' := th + f$ satisfies 
${\rm mon} (h') \subsetneq {\rm mon} (f)$; also by construction $0 \neq h' \in C_r$ and 
${\rm mon} (f) = {\rm mon} (h) \cup {\rm mon } (h')$.  By the induction hypothesis, for any $g \in G$ there are non-zero relations 
$h^g$ and $(h')^g$ in $C_r$ such that ${\rm mon} (h^g) = g \cdot {\rm mon} (h)$ and 
${\rm mon} \bigl( (h')^g \bigr) = g \cdot {\rm mon} (h')$.  Denote by $f^g$ a linear combination of $h^g$ and $(h')^g$ where no 
cancellation of monomials occurs; clearly we have ${\rm mon} (f^g) = g \cdot {\rm mon} (f)$.\end{proof}
\begin{theorem} \label{thm:DPGfan} For a general choice of points $p_1,\dots, p_r$, the Weyl group $W_r$ acts by symmetries on the Gr\"obner fan of $C_r$.
\end{theorem}
\begin{proof} Since the multigraded Hilbert Series of the ideals $C_r$ are independent of the choice of the points, all monomial initial ideals of all $C_r$'s have a common multigraded Hilbert function, and hence they form an antichain of monomial ideals (under inclusion). By Theorem 1.1 in~\cite{DMc} any such antichain is finite. Thus, the degrees of their minimal generators are uniformly bounded by an integer $d$. Using Lemma~\ref{lem:gfan} and Theorem~\ref{thm:monomialita} we obtain the desired conclusion.\end{proof}
\noindent The genericity condition in Theorem~\ref{thm:monomialita} is necessary: the action of the Weyl group is not monomial on special Del Pezzo surfaces. Consider a Del Pezzo $X_6$ with an Eckart point, that is a point $q$ in which three exceptional curves intersect (as in the figure). 
\begin{center}
\scalebox {0.6}{
\includegraphics{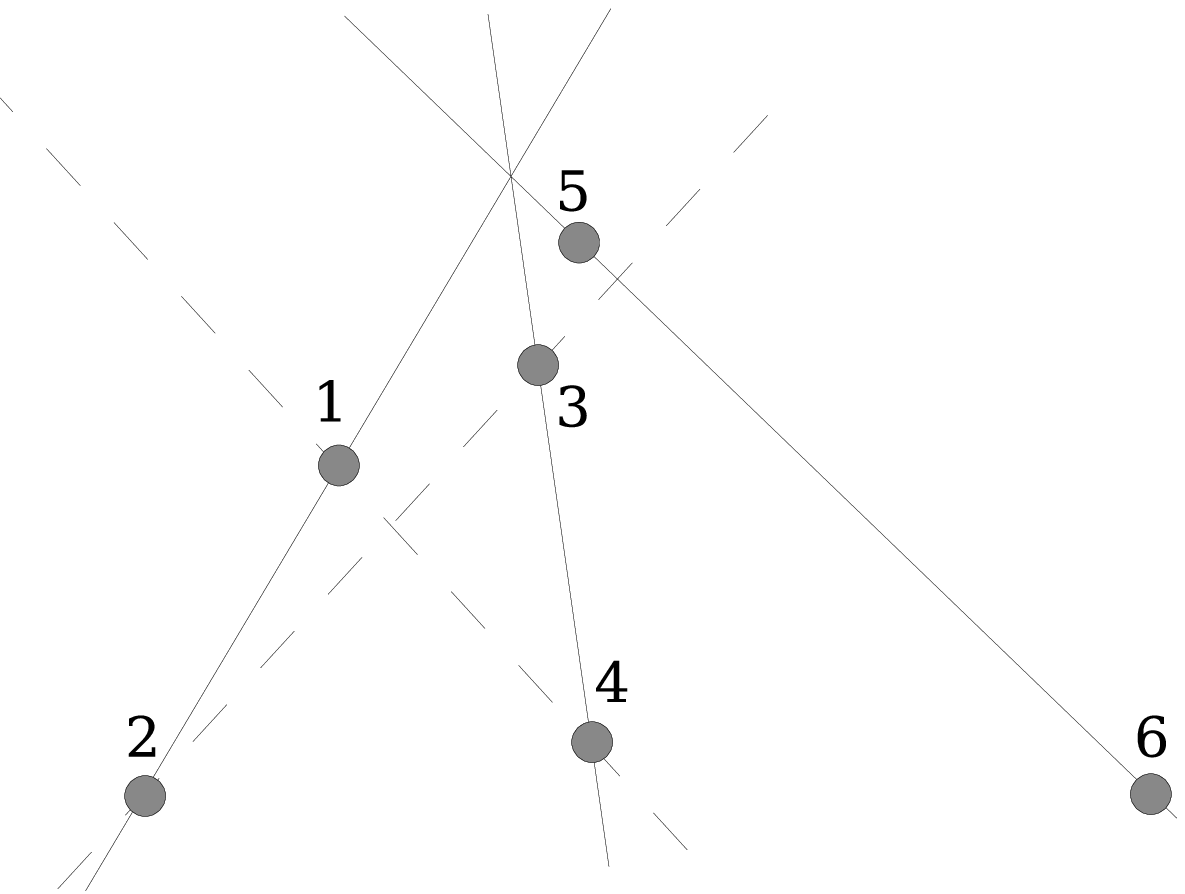}
}
\end{center}
Note that the monomials $f_{12}e_1e_2, f_{34}e_3e_4$ and $f_{56}e_5e_6$ are linearly dependent in $k[E_6]/C_6(p_1,\dots,p_6)$ since so are three lines with a common point in $\mathbb{P}^2$. On the other hand their images under the permutation $(24)\in W_6$:  $f_{14}e_1e_4, f_{23}e_2e_3$ and $f_{56}e_5e_6$ are linearly independent since the corresponding lines do not intersect.   
\section{ Gr\"obner bases and $P$-graded monomial group actions}

In this section the ring of polynomials $R$ is graded by a monoid $P$. The grading by $P$ refines the standard grading by degree and is refined by the grading by monomials. Recall that $G$ acts on $R$ by permuting the variables. Let $\preceq$ be a monomial weight order.
\\
\begin{definition} The action of $G$ is $P$-compatible if for every $g\in G$ the map 
$\hat g:P\rightarrow P$  given by $\hat g(D)=\deg(g(m))$ for any monomial $m$ with 
$\deg(m)=D$ is a well defined endomorphism of $P$.\end{definition}
\noindent
In that case the $\hat g$ are in fact automorphisms of $P$ and the action of $G$ extends to power series via $g(\sum_{D\in P}s_Dt^{D})=\sum_{D\in P}a_Dt^{\hat g(D)}$ and in particular to  Hilbert Series of $P$-homogeneous ideals.
\begin{lemma}Let $I$ be a $P$-homogeneous ideal and assume that $in_{\preceq}(I)$ is generated in degree $\leq d$. If $G$ acts monomially on $I$ up to degree $d$ and the action is $P$-compatible then the Hilbert Series of $I$ is $G$-invariant. \end{lemma}
\begin{proof} By Lemma~\ref{lem:weightsmonomial} we know that 
$g(in_{\preceq}(I))=in_{\preceq_{g^{-1}}}(I)$, thus 
\begin{eqnarray*}
g(HS(I)) & = & \hspace{-1.9pt} g(HS(in_{\preceq}(I))) = HS(g(in_{\preceq}(I))) = 
HS(in_{\preceq_{g^{-1}}}(I)) = \\
& = & \hspace{-1.9pt} HS(I)
\end{eqnarray*}
\end{proof}
\noindent
Note that the $W_r-$action on $k[E_r]$ is ${\rm Pic}(X_r)$ compatible.\\
\\
Now, assume the action of $G$ on $R$ is $P$-compatible and also monomial on $I$. In the remainder of this section we discuss how, under these assumptions, partial knowledge of the Hilbert Series can be used to characterize the initial ideals of $I$. For a vector space $V$ we use $| V |$ to denote its dimension.

\begin{lemma} \label{lem:NZD}
Let $Q\subset J$ be two homogeneous ideals and let $e \in R$ be a nonzerodivisor on $R/Q$ of multidegree $\underline e\in P$. 
Suppose that for some $D \in P$ and $a \in \mathbb{Z}_{\geq 0}$ we have $Q_{D + a \underline e}=J_{D + a \underline e}$; 
then $Q_D=J_D$.
\end{lemma}
\begin{proof}
By induction on $a$ we reduce to the case $a=1$.  Since $Q\subset J$, $Q_D\subset J_D$. Conversely, let $j\in J_D$, 
then $je\in J_{D + \underline e}=Q_{D + \underline e}$ so $j\in Q_D$ since $e$ is a nonzerodivisor on $R/Q$.
\end{proof}

\begin{lemma}\label{lem:HS}
Let $e_1,\dots, e_m$ be variables in $R$ with $\deg(e_i)= \underline e_i$.  Let 
$Q^1,\dots, Q^m$ and $J^1, \dots, J^m$ be homogeneous ideals with $Q^i\subset J^i$ and let 
$s(t)=\sum_{D\in P}S_Dt^{D}$ be a power series. If 
\begin{enumerate}
\item{ All the $J^i$ have the same Hilbert function,  $HS(J^i,t)=s(t)$;}\label{0primo}
\item{ For each $i$, $e_i$ is not a zerodivisor in $R/Q^i$;} \label{0secondo}
\item{ All the $Q^i$ have a common Hilbert function;}
\item{ For all $D\in P$ there exists $j$ and natural numbers $a_i$ such that  
$$|Q^j_{D+\sum_i a_i \underline e_i}|=S_{D+\sum_i a_i \underline e_i};$$}\label{quarto} 
\end{enumerate}
then $Q^i = J^i$ for all $i$.
\end{lemma}

\begin{proof}
Let $D$ be any multidegree in $P$. By (4) there are natural numbers $a_i$ such that, for some $j$, 
\[ |Q^j_{D+\sum_i a_i \underline e_i}|=S_{D+\sum_i a_i \underline e_i} \]
Writing $\sum_i a_i \underline e_i = B + a_j \underline e_j$ with $B = \sum_{i\neq j} a_i \underline e_i$ and using property (\ref{0primo}) 
we conclude
\[ |Q^j_{D + B + a_j \underline e_j}| = S_{D + B + a_j \underline e_j}=|J^j_{D + B + a_j \underline e_j}| \]
Since by property (\ref{0secondo}) $e_j$ is not a zerodivisor in $R/Q^j$, it follows from Lemma~\ref{lem:NZD} and property (\ref{0primo}) that 
\[|Q^j_{D+B}|=|J^j_{D+B}|=S_{D+B}\]
All terms in the above equalities are independent  of $j$ so applying the same argument to the equality
\[|Q^k_{D+B}|=S_{D+B}=|J^k_{D+B}|\]
we can remove from $B$ the $\underline e_k$ component; iterating this process we see that $|Q^l_D|=S_D$ for some index $l$ 
(and hence for all) which proves the statement since $D$ was arbitrary.
\end{proof}

\begin{lemma} \label{lem:ini}
Let $I$ be a homogeneous ideal in $R$ and let $G$ be a group acting monomially on $I$ up to degree $d$ and transitively on 
$e_1,\dots, e_m$. Let $Q \subset {\rm in}_{\preceq}(I)$ be a monomial ideal with generators of total degree at most $d$. If
\begin{enumerate}
\item {The variable $e_1$ is not a zerodivisor in $R/Q$;}
\item {For all $g\in G$, $HS(R/Q,t) = HS(R/g(Q),t)$;}
\item {For all $D\in P$ there are natural numbers $a_i$ such that 
$$ |Q_{D + \sum_i a_i \underline e_i}| = |I_{D + \sum_i a_i \underline e_i}|;$$}
\end{enumerate} 
then $Q = {\rm in}_{\preceq}(I)$.
\end{lemma}

\begin{proof} Choose $g_j\in G$ such that $g_j(e_1)=e_j$. Let $Q^j=g_j(Q)$ and note that, by property (2), $HS (Q^j, t) = HS (Q , t)$ for all $j$.\\
Let $J^j=in_{\preceq_{g_j^{-1}}}(I)$;  since $G$ acts monomially on $I$ and $Q$ has generators of degrees at most $d$ then $Q^j\subset J^j$. Moreover, $HS(J^j , t) = HS (I , t)$, since $J^j$ is an initial ideal of $I$.
Applying Lemma~\ref{lem:HS} with $s(t)=HS(I,t)$ we conclude that $Q = {\rm in}_{\preceq}(I)$.
\end{proof}

We now specialize to $R=k[E_r]$, $G=W_r$, $I=C_r(p_1,\dots,p_r)$ and assume for the rest of the section that the points $p_1,\dots, p_r$ have been chosen general enough so that $W_r$ acts monomially on $C_r$ up to degree $d$ (Theorem~\ref{thm:monomialita}).\\ 
\\
Recall that $Cox(X_r)\cong k[E_r]/C_r$ and that the multigraded Hilbert Series of $Cox(X_r)$ is given by \[HS(Cox(X_r),t)=\sum_{D\in {\rm Pic}(X_r)} h^0({\cal O}_X[D])t^{D}\] 
so in particular it is $W_r$ invariant (since $h^0({\cal O}[D])$ depends only on intersection products which are $W_r$ invariant).\\ 
\noindent
\begin{lemma} \label{lem:sim} Let $M^1\subseteq in_{\preceq}(C_r)$ be a monomial ideal with generators of degree $\leq d$ which does not involve the variable $e_1$; then $M^1=in_{\preceq}(C_r)$ if and only if the following two conditions are satisfied:
\begin{enumerate}
\item{$HS(M^1,t)$ is $W_r$-invariant;}\label{primo}
\item{For $m\in\mathbb{N}$ and all $a_i \gg 0$, 
$\bigl| k[E_r]/M^1 \bigr|_{m\ell+\sum_ia_ie_i}=\binom{m+2}{2}$.}\label{secondo}
\end{enumerate} 
\end{lemma}
\begin{proof} Apply Lemma~\ref{lem:ini} with $Q=M^1$. Hypotheses (1) and (2) of Lemma~\ref{lem:ini} follow immediately from our assumptions thus we verify hypothesis (3). For that let $D$ be any divisor and note that there exist natural numbers $a_i \gg 0$ such that
\[D+\sum a_ie_i=m\ell+\sum b_ie_i\]
with $b_i\geq 0$. Thus $h^0({\cal{O}}[D+\sum{a_ie_i}])=h^0({\cal O}[m\ell]) = \binom{m+2}{2}$ and
by~\ref{secondo} above we have
\[ \bigl| k[E_r]/M^1 \bigr|_{D+\sum a_ie_i}=\binom{m+2}{2}=|Cox(X_r)|_{D+\sum a_ie_i}\]
Thus, for any $D\in {\rm Pic}(X_r)$ there are natural numbers $a_i$ such that
\[|M^1|_{D+\sum a_ie_i}=|C_r|_{D+\sum a_ie_i}\]
as we wanted to prove.
If conversely $M^1=in_{\preceq}(C_r)$,~\ref{primo} and~\ref{secondo} are obvious from the remark preceding this Lemma.
\end{proof}
The third condition in the last Theorem is easy to verify.
\begin{lemma} \label{lem:eq} Let $M$ be a monomial ideal in $k[E_r]$ and let $D = m\ell+a_1e_1+\dots +a_re_r \in P$.  If for $m\in \mathbb{N}$ and $a_i \gg 0$ the dimensions $ \bigl| k[E_r]/M \bigr|_D$ or 
$\bigl| k[E_r]/(M:(e_1\dots e_r)^{\infty}) \bigr|_D$ do not depend on the $a_i$'s, then 
$$ \bigl| k[E_r]/M \bigr|_D = \bigl| k[E_r]/(M:(e_1\dots e_r)^{\infty}) \bigr|_D $$
\end{lemma}
\begin{proof} Let $u=e_1e_2\dots e_r$ and let $k$ be such that $(M:u^k) = (M:u^{\infty})$. Consider the exact sequence
\[0\rightarrow k[E_r]/(M:u^k)[-\deg(u^k)]\rightarrow k[E_r]/M\rightarrow k[E_r]/(M+(u^k))\rightarrow 0\]
and note that $\bigl| k[E_r]/(M+(u^k)) \bigr|_D = 0$ for all $m$ and all $a_i$ sufficiently large. As a result the first map determines an isomorphism between the corresponding graded components.
\end{proof}
\section{The homogeneous coordinate rings for Del Pezzo surfaces}
As an application of the techniques developed so far, we prove in this section that for $r=4,5$ and cubic surfaces without Eckart points the ideal of relations  $C_r(p_1,\dots, p_r)$ of the Cox ring of every Del Pezzo surface $X_r$ is generated by the quadrics $Q_r(p_1,\dots, p_r)$ coming from conic divisor classes.\\
As remarked earlier this fact was conjectured by Batyrev and Popov in~\cite{BP}.\\
\begin{theorem}\label{thm:BPgen}  Assume that $M\subset in_{\preceq}(Q_r(p_1,\dots,p_r))$ is a monomial ideal which does not involve all variables. If $p_1,\dots, p_r$ are sufficiently general then the following are equivalent:
\begin{enumerate}
\item{$M=in_{\preceq}(Q_r)$ and $Q_r=C_r$;} \label{prim}
\item{\begin{enumerate}
\item{$HS(M,t)$ is $W_r$-invariant,}
\item{For $m\in\mathbb{N}$ and all $a_i \gg 0$, 
$\bigl| k[E_r]/M^1 \bigr|_{m\ell+\sum_ia_ie_i}=\binom{m+2}{2}$.}
\end{enumerate}} \label{second}
\end{enumerate}
\end{theorem}
\begin{proof} Clearly~\ref{prim} implies~\ref{second} since these properties are obviously satisfied by the Hilbert Series of the Cox ring. If~\ref{second} holds, note that
\[M\subseteq in_{\preceq}(Q_r)\subseteq in_{\preceq}(C_r)\]
so that by Lemma~\ref{lem:sim}, $M=in_{\preceq}(C_r)$. Thus $Q_r$ and $C_r$ have a common initial ideal and in particular the same Hilbert Series. Since $Q_r\subseteq C_r$ the ideals coincide. Note that the genericity of $p_1,\dots, p_r$ was used to guarantee that $W_r$ acted monomially on $C_r$ up to degree $d$, a necessary condition for Lemma~\ref{lem:sim}.\end{proof}
Let $S\subset (\mathbb{P}^2)^r$ be the open set of $r$-tuples of points in general position (no three on a line, no six on a conic) and let $U\subset S$ be open and nonempty.
\begin{theorem} \label{thm:BP} Assume that $M_r$ is a monomial ideal which does not involve all variables. If 
\begin{enumerate}
\item{$M_r\subset in_{\preceq}(Q_r(p_1,\dots, p_r))$ for all $(p_1,\dots ,p_r)\in U$;} 
\item{$HS(M_r,t)$ is $W_r$-invariant;}
\item{For $m\in\mathbb{N}$ and all $a_i \gg 0$, 
$\bigl| k[E_r]/M^1 \bigr|_{m\ell+\sum_ia_ie_i}=\binom{m+2}{2}$;}
\end{enumerate}
then $C_r=Q_r$ for all Del Pezzo surfaces $X_r$ obtained by blowing up $\mathbb{P}^2$ at points $(p_1,\dots,p_r)\in U$. Moreover $M_r$ is a common initial ideal of all $C_r(p_1,\dots,p_r)$ for $(p_1,\dots, p_r)\in U$.
\end{theorem}
\begin{proof} Since $E$ is irreducible, $U$ and the open set on which the action $W_r$ is monomial guaranteed by Theorem~\ref{thm:monomialita} must intersect. Thus, for a sufficiently general choice of points $p_1,\dots, p_r$, Theorem~\ref{thm:BPgen} implies that $M_r=in_{\preceq}(C_r(p_1,\dots, p_r))$. Since the Multigraded Hilbert Series of the Cox ring is independent of the coordinates of the points we conclude that $M_r=in_{\preceq}(C_r(p_1,\dots, p_r))$ for {\it every} Del Pezzo surface obtained by blowing up $\mathbb{P}^2$ at points in $U$.\\ 
On the other hand the inclusions
\[Q_r\subseteq C_r ~~\text{ and }~~ in_{\preceq}(C_r)=M_r\subseteq in_{\preceq}(Q_r)\]
imply that the initial ideals (and hence the Hilbert Series) of $C_r$ and $Q_r$ coincide for all Del Pezzo surfaces $X_r$ for $(p_1,\dots, p_r)\in U$ and hence that $C_r=Q_r$ for all $(p_1,\dots, p_r)\in U$.
\end{proof}
We now construct the monomial ideals $M_r$.
\subsection{The Del Pezzo surface $X_4$}
In this section $\preceq$ denote the reverse lexicographic on $k[E_4]$ with
\[ e_1 \succ e_2 \succ e_3 \succ e_4 \succ f_{12} \succ f_{13} \succ f_{23} \succ f_{14} \succ f_{24} \succ f_{34}\]
For a finite subset $A\subset {\rm Pic}(X_r)$ and $d\in {\rm Pic}(X_r)$ let $A+d=\{a+d: a\in A\}$, $2A=\{a+a:a\in A\}$ and  $t^A=\sum_{a\in A}t^a$.
\begin{lemma} \label{lem:M4} The ideal $M_4=(f_{23}f_{14}, e_1f_{14}, e_1f_{13}, e_2f_{12}, e_1f_{12})$ has the following properties:\begin{enumerate}
\item{The generators of $M_4$ do not involve $e_4$;}
\item{$M_4\subseteq in_{\preceq}(Q_4)$;}
\item{The Multigraded Hilbert Series of $k[E_4]/M_4$ is $W_4$-invariant;}
\item{For $m\in\mathbb{N}$ and all $a_i \gg 0$, 
$|k[E_r]/M_4|_{m\ell+\sum_ia_ie_i}=\binom{m+2}{2}$.}
\end{enumerate}
\end{lemma}
\begin{proof} (1.) and (2.) are obvious (see generating set of $Q_4$ in Section 4). Direct computation shows that the Hilbert Series of $k[E_4]/M_4$ is given by
\[1-t^{C}+t^{-K-C}-t^{-K}\]
where $C$ denotes the set of conic bundles in $X_4$.  Hence (3.) follows since each of the summands is $W_4$-invariant.\\ 
Finally, Lemma~\ref{lem:eq} shows that it suffices to verify (4.) on 
\begin{eqnarray*}
k[E_4]/ (M_4:e_1\dots e_4) & = & k[E_4] / (f_{12} , f_{13} , f_{14}) = \\
& = & k [e_1 , e_2 , e_3 , e_4 , f_{23} , f_{24} , f_{34}]
\end{eqnarray*}
In this ring, the monomials $e_1^{s_1}e_2^{s_2}e_3^{s_3}e_4^{s_4}f_{23}^{h_{23}}f_{24}^{h_{24}} f_{34}^{h_{34}}$ of multidegree $m \ell + a_1 \underline e_1 + \dots + a_4 \underline e_4$ with $m , a_i \geq 0$ correspond to (non-negative) integral solutions of the system of equations:
\begin{eqnarray*}
s_1 & = & a_1\\[-3pt]
s_2 - h_{23} - h_{24} & = & a_2\\[-3pt]
s_3 - h_{23} - h_{34} & = & a_3\\[-3pt]
s_4 - h_{24} - h_{34} & = & a_4\\[-3pt]
h_{23} + h_{24} + h_{34} & = & m
\end{eqnarray*}
There are $\binom{m + 2} {2}$ such solutions (for all $a_i \gg 0$) since they are determined by decompositions of $m$ as a sum of three non-negative integers.  \end{proof}
\begin{cor} $Cox(X_4)\cong k[E_4]/Q_4$.\end{cor}
\begin{proof} Follows immediately from Lemma~\ref{lem:M4} and Theorem~\ref{thm:BP}.\end{proof}

\subsection{The Del Pezzo surfaces $X_5$}
Let $\preceq$ be any monomial order on $k[E_5]$ refining the one determined by 
the following weights:
\begin{equation} \label{pesi}
\begin{array}{ccccccccccccccccc}
e_4 & f_{13} & e_2 & f_{23} & f_{15} & g & e_1 & f_{25} \\
13 & 11 & 10 & 9 & 8 & 8 & 7 & 7 \\[9pt]
e_3 & e_5 & f_{12} & f_{14} & f_{24} & f_{34} & f_{35} & f_{45} \\
6 & 6 & 6 & 6 & 6 & 6 & 6 & 1
\end{array}
\end{equation}
and let $M_5$ denote the edge ideal of the following graph (the ideal generated by products of variables which are connected by an edge). 
\begin{center}
\scalebox {0.4}{
\includegraphics{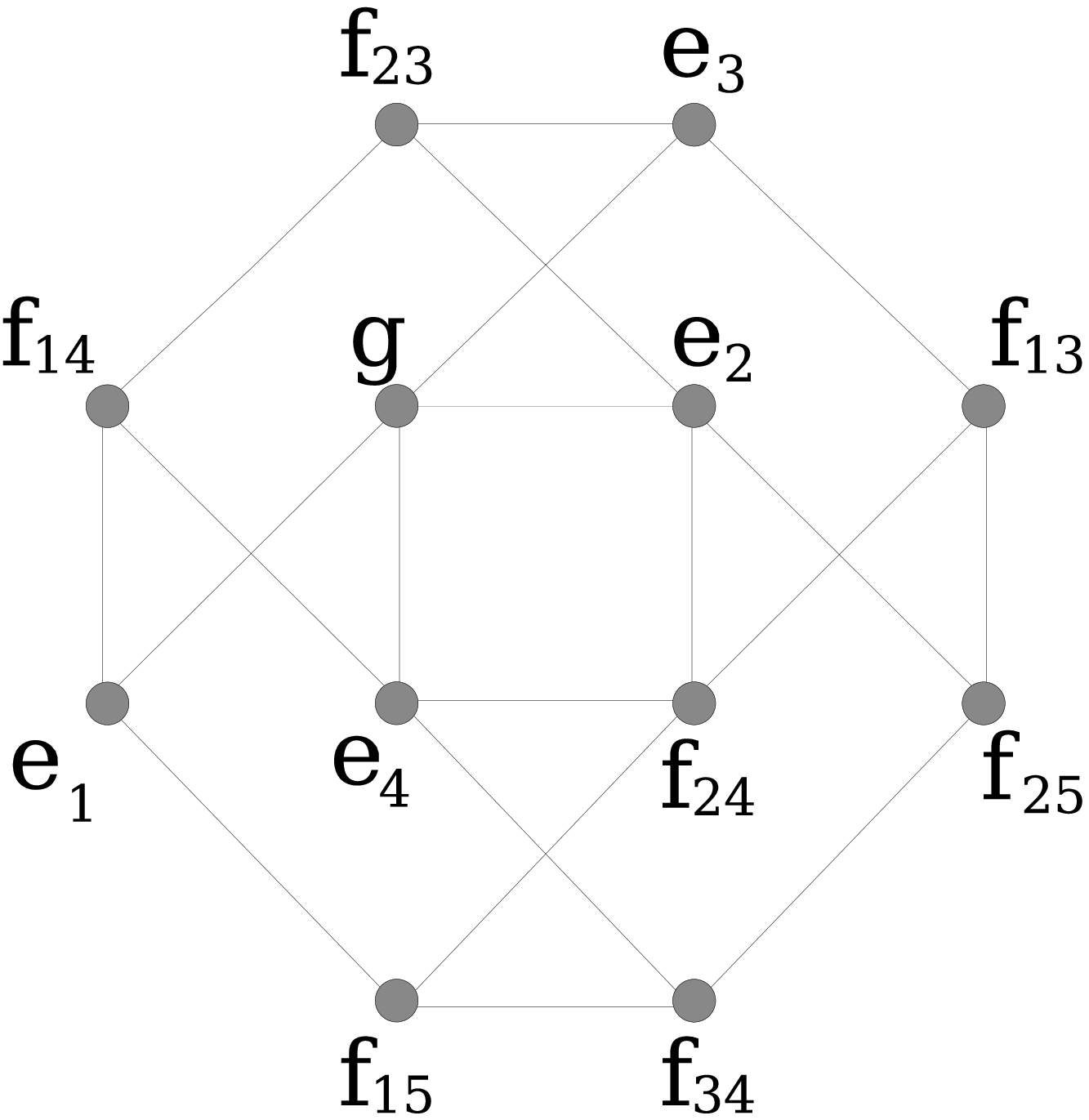}
}
\end{center}
\begin{lemma} \label{lem:M5a}  The ideal $M_5$ is contained in the $\preceq$ initial ideal of $Q_5$ for every Del Pezzo surface $X_5$.
\end{lemma}
\begin{proof} On a Del Pezzo surface any three sections of a conic bundle $D$ are linearly dependent. As a result, given any three monomials of degree $D$ and any choice of points $p_1,\dots,p_5$ some linear combination of the monomials lies in $Q_5(p_1,\dots, p_5)$. Moreover the coefficients of this linear combination must be all nonzero since any two distinct distinguished global sections (i.e. monomials) of $D$ have different support and are therefore independent.\\
As a result, to show that a given monomial $m$ of degree $D$  is a $\preceq$-initial term of $Q_5(p_1,\dots,p_5)$ for all Del Pezzo surfaces $X_5$ it suffices to find two more monomials of degree $D$ which are $\preceq$-smaller.\\    
The $\preceq$ weights of all monomials in each conic $D$ are recorded in the table. 
$$ \begin{array} {|c|cc|cc|cc|cc|}
\hline 
{\rm Pic} (X)-{\rm degree} & \multicolumn {8}{c|}{\text{Monomials and their weights}} \\ \hline 
\ell - e_1 & e_4 f_{14} & 19 & e_3 f_{13}     & 17 & e_2 f_{12}    & 16 & e_5 f_{15} & 14 \\[2.5pt]
\ell - e_2 & e_4 f_{24} & 19 & e_3 f_{23}     & 15 & e_1 f_{12}    & 13 & e_5 f_{25} & 13 \\[2.5pt]
\ell - e_3 & e_4 f_{34} & 19 & e_2 f_{23} & 19 & e_1 f_{13} & 18 & e_5 f_{35} & 12 \\[2.5pt]
\ell - e_4 & e_2 f_{24} & 16 & e_1 f_{14}     & 13 & e_3 f_{34}    & 12 & e_5 f_{45} &  7 \\[2.5pt]
\ell - e_5 & e_2 f_{25} & 17 & e_1 f_{15} & 15 & e_4 f_{45} & 14 & e_3 f_{35} & 12 \\[2.5pt]
2 \ell - e_1 - e_2 - e_3 - e_4 & f_{24} f_{13} & 17 & f_{14} f_{23}     & 15 & e_5 g    & 14 & f_{12} f_{34} & 12 \\[2.5pt]
2 \ell - e_1 - e_2 - e_3 - e_5 & e_4 g & 21 & f_{13} f_{25} & 18 & f_{23} f_{15} & 17 & f_{12} f_{35} & 12 \\[2.5pt]
2 \ell - e_1 - e_2 - e_4 - e_5 & e_3 g & 14 & f_{24} f_{15} & 14 & f_{14} f_{25}    & 13 & f_{12} f_{45} &  7 \\[2.5pt]
2 \ell - e_1 - e_3 - e_4 - e_5 & e_2 g & 18 & f_{15} f_{34} & 14 & f_{13} f_{45} & 12 & f_{14} f_{35} & 12 \\[2.5pt]
2 \ell - e_2 - e_3 - e_4 - e_5 & e_1 g & 15 & f_{25} f_{34} & 13 & f_{24} f_{35} & 12 & f_{23} f_{45} & 10 \\
\hline 
\end{array} $$
Note that the generators of $M_5$ are the first two monomials in each row. They have higher weights than the last two and hence $M_5\subset in_{\preceq}(Q_5(p_1,\dots, p_5))$ for all Del Pezzo surfaces $X_5$.
\end{proof}
\begin{lemma} \label{lem:M5b} The ideal $M_5$ has the following properties: 
\begin{enumerate}
\item{The generators of $M_5$ do not involve the variable $e_5$;}
\item{$HS(M_5,t)$ is $W_5$-invariant;}
\item{For $m\in\mathbb{N}$ and all $a_i \gg 0$, 
$\bigl| k[E_5]/M_5 \bigr|_{m\ell+\sum_ia_ie_i}=\binom{m+2}{2}$.}
\end{enumerate}
\end{lemma}
\begin{proof} To verify (2.) we describe the multigraded Hilbert Series of $k[E_5]/M_5$ explicitly. Consider the following sets of divisors on $X_5$ (note that each set is invariant under the action of $W_5$).\\
\\
\bigskip
\begin{tabular}{c|c}
Set & Divisors\\
\hline
\hline
$C$ & Conic bundles\\
$D$ & Divisors $d=c+v$ for $c\in C$, $v\in E_5$ with $c\cdot v=1$\\ 
$F$  & Divisors $f=2c+v$ for $c\in C$, $v\in E_5$ with $c\cdot v=1$\\
$G$ & Divisors $g=c_1+c_2$ for $c_i\in C$, $c_1\cdot c_2=2$\\
$H$ & Divisors $h=2c_1+c_2$ for $c_i\in C$, $c_1\cdot c_2=2$\\
$J$  & Divisors $j=c_1+c_2+c_3$ for $c_i\in C$, $c_i\cdot c_j=2$ for $i\neq j$\\  
$K$ & Canonical divisor
\end{tabular}\\
\\
Now let  \[\alpha=1-2t^C+3t^D+t^{2C}-t^F-3t^{-K}-6t^{G}\]
Direct computer computation shows that the numerator of the Hilbert Series of $k[E_5]/M_5$ is given by:
\begin{equation} \label{sehil}
(\alpha + \alpha^{*}) + t^{J} + 12t^{H}
\end{equation}
where $(-)^*$ is linear and $(t^A)^*$ denotes $t^{-3K-A}$.\\
It follows that the Hilbert Series is invariant under the action of $W_5$ since each term of (\ref{sehil}) is.\\
By Lemma~\ref{lem:eq} we can verify (3.) in the quotient
\begin{eqnarray*}
& k[E_5]/\bigl(M_5:(e_1\dots e_5)^{\infty}\bigr)=\\
& k[E_5] / \bigl( f_{13} , f_{23} , f_{14} , f_{24} , f_{34} , f_{15} , f_{25} , g \bigr)=\\ 
& =k [e_1, e_2, e_3, e_4, e_5 , f_{12} , f_{35} , f_{45} ]
\end{eqnarray*}
The monomials $e_1^{s_1} e_2^{s_2} e_3^{s_3} e_4^{s_4} e_5^{s_5} f_{12}^{h_{12}} 
f_{35}^{h_{35}} f_{45}^{h_{45}}$ in this ring of multidegree 
$m \ell + a_1e_1 + \dots +a_5e_5$ with $m , a_i \geq 0$ correspond to non-negative integral solutions of the system of equations:
$$ \begin{array}{ccc}
s_1 - h_{12} & = & a_1\\
s_2 - h_{12} & = & a_2\\
s_3 - h_{35} & = & a_3\\
s_4 - h_{45} & = & a_4\\
s_5 - h_{35} - h_{45} & = & a_5\\
h_{12} + h_{35} + h_{45} & = & m
\end{array} $$
For all sufficiently large $a_i$ these solutions correspond to ways of writing $a$ as a sum of three natural numbers: there are $\binom{m + 2} {2}$ such monomials.\\
\end{proof}
\begin{cor} For all Del Pezzo surfaces $X_5$ of degree 4 we have 
$C_5(p_1,\dots, p_5)=Q_5(p_1,\dots, p_5)$. In particular
\[Cox\bigl(X_5(p_1,\dots,p_5)\bigr)\cong k[E_5]/Q_5(p_1,\dots, p_5)\]\end{cor}
\begin{cor} In the notation of Lemma~\ref{lem:M5b} the multigraded Hilbert Series of $Cox(X_5)$ is given by 
\[HS(Cox(X_5))=\alpha+t^{J}+12t^{H}+\alpha^{*}\]
\end{cor}
\begin{cor} The ideal $M_5$ is an initial ideal of $C_5(p_1,\dots, p_5)$ for all Del Pezzo surfaces $X_5$. In particular the rings $Cox(X_r)$ for $r\leq 5$ are Koszul algebras.\end{cor}
\begin{proof} The corollaries follow immediately from Lemmas~\ref{lem:M5a} and~\ref{lem:M5b} and Theorem~\ref{thm:BP}. The last statement follows from the well known fact that every $G$-quadratic algebra is Koszul.\end{proof} 

\subsection{The Del Pezzo surfaces $X_6$}
Recall that every nef divisor class $D$ on a Del Pezzo surface can be written as
\[D=n_8 (-K_{X_8}) + \ldots + n_2 (-K_{X_2}) + D' \]
where $X_8 \rightarrow X_7 \rightarrow \ldots \rightarrow X_1$ is a sequence of contractions of
$(-1)-$cur\-ves, $K_{X_i}$ is the canonical divisor of $X_i$, the $n_i\geq 0$ and $D'$ is a nef
divisor on $X_1$, a surface which is either the blow up of $\mathbb{P}^2$ at one point or $\mathbb{P}^1\times \mathbb{P}^1$.\\
Using this fact we show that, for every Del Pezzo surface $X_6$, the ideals $Q_6$ and $C_6$ coincide in those degrees $D$ with $D^2=1$ and $-K\cdot D=3$. 
\begin{lemma} \label{lem:gradotre}
Let $X_6$ be a Del Pezzo surface of degree 3 and let $D \in {\rm Pic} (X_6)$ be a nef divisor such that $D^2 = 1$ and $-K \cdot D = 3$.  
Then there is a morphism $\pi : X_6 \rightarrow \mathbb {P}^2$ exhibiting $X_6$ as the blow up of $\mathbb {P}^2$ at 
$6$ points such that $D = \pi ^* \ell $, where $\ell $ is the divisor class of a line in $\mathbb {P}^2$.
\end{lemma}
\begin{proof}
Write 
$$ D = n_6 (-K_{X_6}) + \ldots + n_2 (-K_{X_2}) + D' $$
with $n_j \geq 0$ and $D'$ a nef divisor on either the blow up of $\mathbb{P}^2$ at one point or on $\mathbb{P}^1\times \mathbb{P}^1$.  Since all 
the divisors appearing on the right are nef, all intersection numbers among them are non-negative.  In particular 
the condition $D^2 = 1$ implies $n_6 = \ldots = n_2 = 0$.\\ 
Now, $\mathbb{P}^1\times\mathbb{P}^1$ admits no divisor of square 1 (its Picard group is generated by the rulings $h_1, h_2$ with $h_1^2=h_2^2=0$ and $h_ih_j=1$) and the blow up of $\mathbb{P}^2$ at one point has $D=\ell$ as only solution of the equations $D^2=1$ and $-K\cdot D=3$.
Thus $D$ is the pull-back of the divisor class of a line under a birational morphism 
$X \rightarrow \mathbb {P}^2$.
\end{proof}

\begin{lemma} \label{lem:relations}
Let $X$ be a Del Pezzo surface of degree at most 6.  The ideals $Q_6$ and $C_6$ coincide in all 
degrees $D \in P$, where $D$ is a nef divisor satisfying $D^2 = 1$ and $K \cdot D = -3$.
\end{lemma}

\begin{proof}
The Riemann-Roch formula and the Kodaira vanishing Theorem imply that $\dim Cox (X)_D = 3$.  
Since $k[E_r]/Q_r \rightarrow k[E_r]/C_r = Cox (X)$ is surjective, the result follows if we show 
$\dim \bigl( k[E_r]/Q_r \bigr)_D \leq 3$.

Lemma~\ref{lem:gradotre} implies that we can find a morphism $\pi : X \rightarrow \mathbb {P}^2$ 
which is the blow up of $m$ points $p_1, \ldots , p_m$, such that $D$ is the pull-back of the 
divisor class of a line in $\mathbb {P}^2$.  Denote by $\underline g_1 , \dots , \underline g_m$ 
be the divisor classes of the exceptional divisors of $\pi $ and denote by $g_i$ the variable of 
$k[E_r]$ corresponding to $\underline g_i$ and by $h_{ij}$ the variable corresponding to the strict 
transform of the line through the two blown up points $p_i$ and $p_j$.  Any element $f \in k[E_r]$ 
of degree $D$ is a linear combination of monomials $h_{ij} g_i g_j$, with $1 \leq i < j \leq m$.  
It is enough to show that $f$ is a linear combination of $h_{ij} g_i g_j$, with 
$1 \leq i < j \leq 3$, modulo $Q_r$.

Since $Q_r$ and $C_r$ agree up to coarse degree 2, the ideal $Q_r$ contains relations of the form 
$q = q_{ij} h_{ij} g_j + q_{ri} h_{ri} g_r + q_{si} h_{si} g_s$, where 
$q_{ij} , q_{ri} , q_{si}$ are nonzero constants.  Thus if $i < j$ and $j \geq 4$, we can choose distinct 
$r,s \in \{1,2,3\} \setminus \{i\}$ and use the relation $qg_i$ and induction on $j$ to achieve 
our goal.
\end{proof}
Let $\preceq$ be the monomial order obtained by refining the following weight vector with the reverse lexicographic order:
\[\begin{array}{|ccccccccc|}
\hline
3 & 3 & 3 & 2 & 2 & 2 & 2 & 2 & 2\\ 
f_{23} & g_3 & e_2 & g_2 & g_6 & e_6 & f_{26} & f_{24} & g_4\\
\hline
\hline
2 & 2 & 2 & 2 & 2 & 2 & 2 & 2 & 2\\ 
f_{45} & e_5 & f_{34} & f_16 & e_4 & f_{14} & f_{12} & g_1 & f_{46}\\
\hline
\hline
2 & 2 & 2 & 2 & 2 & 2 & 1 & 1 & 1\\ 
e_3 & f_{36} & f_{35} & f_{13} & f_{56} & f_{25} & e_1 & g_5 & f_{15}\\
\hline
\end{array}\]
We now construct a monomial ideal $M_6$ which is a $\preceq$-initial ideal of $C_r$ for every cubic surface without Eckart points. The tables mentioned in its description are contained in the Appendix (Section~\ref{sec:appendix}).\\  
Let $M_6$ be the monomial ideal generated by the 81 quadratic monomials in the first three columns of Table T1, the 34 cubic monomials in the first column of Table T2 and the first cubic monomial in Table T3. Note that the monomials of each Picard degree $D$ have been written in decreasing $\preceq$-order in the Tables.\\
\\
Recall that an Eckart point on a Del Pezzo surface $X_6$ is a point in which 3 exceptional curves intersect.
\begin{lemma} \label{lem:M6a} The inclusion $M_6\subset in(Q_6(p_1,\dots,p_6))$ holds for every Del Pezzo surface $X_6$ without Eckart points.\end{lemma}
\begin{proof}  The quadratic generators of $M_6$ are initial terms of $Q_6$ for every Del Pezzo surface $X_6$ (appearing in the relation which involves them and the two last monomials in each row of Table T1 in the Appendix).\\ 
For the cubic generators of $M_6$ there are two cases depending on their multidegree $D$.\\
If $D^2=1$ then $h^0(O[D])=3$ so every four monomials in $D$ are linearly dependent.
Moreover, all coefficients of any linear dependence relation are nonzero since if $X_6$ has no Eckart points, Lemma~\ref{lem:gradotre} implies that every three monomials in $D$ are linearly independent. 
In particular (by Lemma~\ref{lem:relations}) the first term in each row of Table T2 is an element of $in_{\preceq}(Q_6(p_1,\dots,p_6))$.\\
Finally, it follows from Lemma~\ref{lem:cubicgenerator} (in the Appendix) that the monomial
$f_{46}f_{13}f_{25}$, the only cubic generator of $M_6$ of degree $-K$, lies in $in_{\preceq}(Q_6)$ for every cubic surface without Eckart points.
\end{proof}
Note that a general Del Pezzo surface $X_6$ has no Eckart points.\\
\\
Direct calculations using Macaulay2 show that:
\begin{lemma} \label{lem:M6b} The ideal $M_6$ has the following properties: 
\begin{enumerate}
\item{The generators of $M_6$ do not involve the variable $g_5$;}
\item{$HS(M_6,t)$ is $W_6$-invariant;}
\item{For $m\in\mathbb{N}$ and all $a_i \gg 0$,
$|k[E_6]/M_6|_{m\ell+\sum_ia_ie_i}=\binom{m+2}{2}$.}
\end{enumerate}
\end{lemma}
\begin{cor} For all Del Pezzo surfaces $X_6$ without Eckart points, $C_6=Q_6$. In particular
\[Cox\bigl(X_6(p_1,\dots,p_5)\bigr)\cong k[E_6]/Q_6(p_1,\dots, p_6)\]\end{cor}
\begin{proof} Follows from Lemmas~\ref{lem:M6a} and~\ref{lem:M6b} and Theorem~\ref{thm:BP} where $U$ is the set of points $(p_1,\dots, p_r)$ in general position such that the corresponding blow up of $\mathbb{P}^2$ is a Del Pezzo surface without Eckart points.\end{proof}
\section{Quadratic Gr\"obner basis}
For $r=4$ and $5$, the quadratic initial ideals which we have exhibited are edge ideals of subgraphs of the graphs of exceptional curves $L_r$. As we will show, this is no accident.\\
\\
Note that the edges of $L_r$ can be colored by the conic bundles by assigning $\deg(v_1v_2)$ to the edge $\{v_1,v_2\}$. Each color class contains exactly $r-1$ edges.
\begin{lemma} \label{lem:QGB} If $M$ is a quadratic initial ideal of $C_r$ then $M$ is the edge ideal of  a subgraph of the graph of exceptional curves with $r-3$ edges on each color class.
\end{lemma}
\begin{proof} If $D=\ell-e_1$, a basis for $R_D$ is given by the $r-1$ monomials $f_{1i}e_i$, for $i=2,\dots, r$. Since $W_r$ acts transitively on conic bundles $D$ we see that $|R_D|=r-1$, and that the monomials in this graded component correspond to pairs of (-1)-curves which intersect, that is, to edges of the graph $H_r$ on the conic bundle $D$.\\  
Now, $|R/M|_D=|R/C_r|_D=h^0({\cal O}[D])=2$ so $|M|_D=r-3$ for all conic bundles. 
If $-K\cdot D=2$ and $D$ is effective but not a conic bundle then $D=C_1+C_2$ and the curves $C_1,C_2$ either do not intersect or are equal. As a result $C_j\cdot D<0$ and $R_D$ is spanned by the monomial $c_1c_2$. Since $h^0({\cal O}[D])=1$, $|M_D|=0$ and the statement follows.
\end{proof} 
In view of Lemma~\ref{lem:QGB} and Theorem~\ref{thm:BP} it becomes a question of interest to characterize all subgraphs of the graphs of exceptional curves whose multigraded Hilbert Series is $W_r$-invariant. For $r\leq5$ this can be accomplished by direct computer exploration and one can show,
\begin{lemma} There are exactly 18 quadratic initial ideals of $C_5$ up to the action of $W_5=D_5$.\end{lemma}
\noindent The corresponding weight vectors are in Table T4 in the Appendix.\\ 
\\
For $r=6$ the space of possible subgraphs is much larger and exhaustive exploration is simply unfeasible. One can in fact show that there are no ``small'' quadratic initial ideals, that is, initial ideals whose generators involve less than 25 variables. This observation depends on the geometry of the 27 lines on the cubic and we will prove it.
\begin{lemma} \label{lem:config} Let $A=\{a_1,a_2,a_3\}$ be a set of three distinct exceptional curves. Then the exceptional curves in $A$ do not form a triangle if and only if there exist sets of $(-1)$-curves $C=\{c_1,\dots,c_6\}$ and $H=\{h_2,\dots,h_6\}$ disjoint from $A$ such that
\begin{enumerate}
\item{The $h_i$ are pairwise disjoint;}
\item{$h_i$ intersects $c_i$ and $c_1$ and no other curve in $C$;}
\item{There is a conic bundle $D$ whose distinguished global sections are precisely the $c_ih_i$.}
\end{enumerate}
\end{lemma}
\begin{proof} If the curves in $A$ form a triangle, it is easy to see that every conic bundle $D$ contains exactly one monomial divisible by a variable in $A$ so (3.) is impossible.\\
Conversely note that $C=\{e_1,\dots, e_6\}$ and $H=\{f_{12},f_{13},\dots f_{16}\}$ satisfy conditions $(1.)$, $(2.)$ and $(3.)$ (with $D=\ell-e_1$). Moreover, the sets of curves  $\{f_{23},f_{24},f_{25}\}$, $\{f_{23},f_{46},f_{45}\}$ and $\{f_{23},f_{46},g_1\}$ induce subgraphs with all isomorphism types of graphs of size 3 except the triangle. If the curves in $A$ form any such graph, the action of the Weyl group can carry them into one of these. As a result their complement contains the required sets of exceptional curves.
\end{proof}
\begin{theorem} If there is a quadratic initial ideal $N$ for $Q_6$, its generators must involve at least 25 of the 27 variables.\end{theorem} 
\begin{proof} If $N$ involves 24 variables or less pick three $\{v_1,v_2,v_3\}$ which do not appear in the generators and denote by $S$ the subgraph of $G_6$ that they span.\\
If this subgraph is not a triangle, Lemma~\ref{lem:config} shows that there are sets of curves $C=\{c_1,\dots,c_6\}$ and $H=\{h_1,\dots,h_5\}$ disjoint from $S$ such that $h_i$ intersects $c_i$ and $c_1$ and no other curve in $C$ and the $c_ih_i$ are all  distinguished sections of the conic bundle $D$.\\
By Lemma~\ref{lem:QGB} two distinguished sections of $D$, say $c_1h_1$ and $c_2h_2$ 
are not edges of the subgraph of $L_6$ corresponding to $N$. As a result 
\[c_1 \, , \, c_2 \, , \, c_3 \, , \, c_4 \, , \, c_5 \, , \, h_1 \, , \, h_2 \, , \, v_1 \, , \, v_2 \, , \, v_3\]
is a set of independent vertices (i.e. no two joined by an edge in $G_N$) so ${\rm codim}(N) \leq 17 < {\rm codim}(Q_6)$ and $N$ cannot be an initial ideal.\\
If the subgraph $S$ is a triangle, then the situation is very different and there are subgraphs of $G_6$ with the same Hilbert function as $Q_6$. Computer calculations done with Macaulay2 show that they are not initial ideals for any weight vector. \end{proof}


\section{Appendix: Generators of $M_6$}\label{sec:appendix}
This section contains explicit calculations that were used in the proof of Lemma~\ref{lem:M6a}. All
polynomials have been written in decreasing $\preceq$-order.\\
$M_6$ is generated by the 81 quadratic monomials in the first two columns of Table T1, the 34 cubic monomials in the first column of Table T2 and the first monomial in Table T3.\\
\[
\begin{array}{|c|ccc|cc|}
\hline
\text{\large{T1}} & & & & &\\
\hline
{2\ell-e_1-e_2-e_3-e_5}  &
e_6g_4 & g_6e_4 & f_ { 1 2 }f_ { 3 5 } & f_ { 1 3 }f_ { 2 5 } & f_ { 2 3 }f_ { 1 5 }\\ [1.5pt]
{3\ell  -e_1  -2e_2  -e_3  -e_4  -e_5  -e_6} &
f_ { 2 3 }g_3 & g_6f_ { 2 6 } & f_ { 2 4 }g_4 & f_ { 1 2 }g_1 & f_ { 2 5 }g_5\\ [1.5pt]
{2\ell -e_1  -e_2  -e_3 -e_6}&
f_ { 2 3 }f_ { 1 6 } & g_4e_5 & f_ { 1 2 }f_ { 3 6 } & f_ { 2 6 }f_ { 1 3 } & e_4g_5\\ [1.5pt]
{2\ell  -e_1  -e_2  -e_4  -e_5}&
g_3e_6 & f_ { 4 5 }f_ { 1 2 } & g_6e_3 & f_ { 1 4 }f_ { 2 5 } & f_ { 2 4 }f_ { 1 5 }\\ [1.5pt]
{3\ell  -e_1  -e_2  -2e_3  -e_4  -e_5  -e_6}&
f_ { 2 3 }g_2 & g_4f_ { 3 4 } & g_6f_ { 3 6 } & g_1f_ { 1 3 } & f_ { 3 5 }g_5\\ [1.5pt]
{\ell  -e_1}&
e_2f_ { 1 2 } & e_6f_ { 1 6 } & e_4f_ { 1 4 } & e_3f_ { 1 3 } & e_5f_ { 1 5 }\\ [1.5pt]
{2\ell  -e_1  -e_2-e_4-e_6} &
g_3e_5 & f_ { 2 4 }f_ { 1 6 } & f_ { 2 6 }f_ { 1 4 } & f_ { 1 2 }f_ { 4 6 } & e_3g_5\\ [1.5pt]
{2\ell-e_1-e_3  -e_4  -e_5} &
e_2g_6 & g_2e_6 & f_ { 1 4 }f_ { 3 5 } & f_ { 4 5 }f_ { 1 3 } & f_ { 3 4 }f_ { 1 5 }\\ [1.5pt]
{3\ell -e_1-e_2-e_3 -2e_4 -e_5 -e_6}&
g_3f_ { 3 4 } & g_2f_ { 2 4 } & f_ { 1 4 }g_1 & g_6f_ { 4 6 } & f_ { 4 5 }g_5\\ [1.5pt]
{\ell-e_2} &
f_ { 2 3 }e_3 & e_6f_ { 2 6 } & f_ { 2 4 }e_4 & e_5f_ { 2 5 } & f_ { 1 2 }e_1\\ [1.5pt]
{2\ell  -e_1 -e_2 -e_5 -e_6}&
g_3e_4 & g_4e_3 & f_ { 1 2 }f_ { 5 6 } & f_ { 1 6 }f_ { 2 5 } & f_ { 2 6 }f_ { 1 5 }\\ [1.5pt]
{2\ell  -e_1-e_3-e_4-e_6}&
g_2e_5 & f_ { 3 4 }f_ { 1 6 } & f_ { 1 4 }f_ { 3 6 } & f_ { 4 6 }f_ { 1 3 } & e_2g_5\\ [1.5pt]
{2\ell-e_2  -e_3  -e_4  -e_5}&
f_ { 2 3 }f_ { 4 5 } & e_6g_1 & f_ { 2 4 }f_ { 3 5 } & f_ { 3 4 }f_ { 2 5 } & g_6e_1\\ [1.5pt]
{3\ell  -e_1  -e_2  -e_3  -e_4  -2e_5  -e_6}&
g_3f_ { 3 5 } & g_4f_ { 4 5 } & g_6f_ { 5 6 } & g_2f_ { 2 5 } & g_1f_ { 1 5 }\\ [1.5pt]
{\ell-e_3}&
f_ { 2 3 }e_2 & f_ { 3 4 }e_4 & e_6f_ { 3 6 } & e_5f_ { 3 5 } & f_ { 1 3 }e_1\\ [1.5pt]
{2\ell  -e_1-e_3-e_5-e_6}&
e_2g_4 & g_2e_4 & f_ { 1 6 }f_ { 3 5 } & f_ { 1 3 }f_ { 5 6 } & f_ { 3 6 }f_ { 1 5 }\\ [1.5pt]
{2\ell -e_2 -e_3 -e_4-e_6}&
f_ { 2 3 }f_ { 4 6 } & f_ { 2 6 }f_ { 3 4 } & e_5g_1 & f_ { 2 4 }f_ { 3 6 } & e_1g_5\\ [1.5pt]
{3\ell-e_1-e_2-e_3-e_4-e_5-2e_6}&
g_3f_ { 3 6 } & g_2f_ { 2 6 } & f_ { 1 6 }g_1 & g_4f_ { 4 6 } & f_ { 5 6 }g_5\\ [1.5pt]
{\ell-e_4}&
e_2f_ { 2 4 } & f_ { 4 5 }e_5 & e_6f_ { 4 6 } & f_ { 3 4 }e_3 & f_ { 1 4 }e_1\\ [1.5pt]
{2\ell  -e_2  -e_3 -e_5  -e_6}&
f_ { 2 3 }f_ { 5 6 } & e_4g_1 & f_ { 2 6 }f_ { 3 5 } & f_ { 3 6 }f_ { 2 5 } & g_4e_1\\ [1.5pt]
{2\ell  -e_1-e_4  -e_5-e_6}&
g_3e_2 & f_ { 4 5 }f_ { 1 6 } & g_2e_3 & f_ { 1 4 }f_ { 5 6 } & f_ { 4 6 }f_ { 1 5 }\\ [1.5pt]
{\ell -e_5}&
e_2f_ { 2 5 } & f_ { 4 5 }e_4 & e_3f_ { 3 5 } & e_6f_ { 5 6 } & e_1f_ { 1 5 }\\ [1.5pt]
{2\ell-e_2-e_4 -e_5 -e_6}&
f_ { 2 6 }f_ { 4 5 } & g_1e_3 & f_ { 2 4 }f_ { 5 6 } & f_ { 4 6 }f_ { 2 5 } & g_3e_1\\ [1.5pt]
{\ell-e_6}&
e_2f_ { 2 6 } & e_4f_ { 4 6 } & e_3f_ { 3 6 } & e_5f_ { 5 6 } & f_ { 1 6 }e_1\\ [1.5pt]
{2\ell  -e_1  -e_2  -e_3  -e_4}&
f_ { 2 3 }f_ { 1 4 } & g_6e_5 & f_ { 3 4 }f_ { 1 2 } & f_ { 2 4 }f_ { 1 3 } & e_6g_5\\ [1.5pt]
{2\ell-e_3  -e_4 -e_5 -e_6}&
e_2g_1 & f_ { 4 5 }f_ { 3 6 } & f_ { 4 6 }f_ { 3 5 } & f_ { 3 4 }f_ { 5 6 } & g_2e_1\\ [1.5pt]
{3\ell  -2e_1  -e_2  -e_3  -e_4  -e_5  -e_6}&
g_3f_ { 1 3 } & g_6f_ { 1 6 } & g_4f_ { 1 4 } & g_2f_ { 1 2 } & g_5f_ { 1 5 }\\ [1.5pt]
\hline
\end{array}\]

\[\begin{array} {|c|c|ccc|}
\hline
\text{\large{T2}} & & & &\\
\hline 
3\ell -e_1-e_2-2e_3-e_5-e_6 & f_{36}f_{13}f_{25} & f_{23}f_{36}f_{15} & g_4f_{13}e_1 & f_{35}g_5e_4\\ [1.8pt]
4\ell -e_1-2e_2-2e_3-e_4-2e_5-e_6 & g_1f_{13}f_{25} & f_{23}g_1f_{25} & g_6g_4e_1 & f_{35}f_{25}g_5\\ [1.8pt]
3\ell -e_1-e_2-2e_3-e_4-e_5 & f_{34}f_{13}f_{25} & f_{23}f_{34}f_{25} & g_6f_{13}e_1 & e_6f_{35}g_{5} \\ [1.8pt]
2\ell-e_1-e_2-e_3 & e_5f_{13}f_{25} & f_{23}e_5f_{15} & f_{12}f_{13}e_1 & e_6e_4g_5\\ [1.8pt]
3\ell -e_1-2e_2-e_3-e_4-e_5 & f_{24}f_{13}f_{25} & f_{23}f_{24}f_{15} & g_6f_{12}e_1 & e_6f_{25}g_5 \\ [1.8pt]
3\ell -e_1-2e_2-e_3-e_5-e_6 & f_{26}f_{13}f_{25} & f_{23}f_{26}f_{15} & g_4f_{12}e_1 & e_4f_{25}g_5\\ [1.8pt]
3\ell -e_1-e_2-e_3-e_5-2e_6 & f_{16}f_{36}f_{25} & g_4f_{16}e_1 & e_4f_{56}g_5 & f_{26}f_{36}f_{15}\\ [1.8pt]
4\ell-e_1-e_2-2e_3-e_4-2e_5-2e_6 & g_2f_{36}f_{25} & g_2g_4e_1 & f_{35}f_{56}g_5 & g_1f_{36}f_{15}\\ [1.8pt]
3\ell-e_1-2e_2-e_4-e_5-e_6 &  f_{12}f_{46}f_{25} & g_3f_{12}e_1 & e_3f_{25}g_5 & f_{26}f_{24}f_{15}\\ [1.8pt]
3\ell-e_1-e_2-2e_4-e_5-e_6 &  f_{14}f_{46}f_{25} & g_3f_{14}e_1 & f_{45}e_3g_5 & f_{24}f_{46}f_{15}\\ [1.8pt]
3\ell -e_1-e_2-e_4-e_5-2e_6 & f_{16}f_{46}f_{25} & g_3f_{16}e_1 & e_3f_{56}g_5 & f_{26}f_{46}f_{15}\\ [1.8pt]
4\ell -e_1-2e_2-e_3-e_4-2e_5-2e_6 & g_4f_{46}f_{25} & g_3g_4e_1 & f_{56}f_{25}g_5 & f_{26}g_1f_{15}\\ [1.8pt]
4\ell -e_1-2e_2-e_3-2e_4-2e_5-e_6 & g_6f_{46}f_{25} & g_3g_6e_1 & f_{45}f_{25}g_5 & f_{24}g_1f_{15}\\ [1.8pt]
4\ell -e_1-e_2-e_3-2e_4-2e_5-2e_6 &g_2f_{46}f_{25} & g_3g_2e_1 & f_{45}f_{56}g_5 & g_1f_{46}f_{15}\\ [1.8pt]
3\ell -e_1-e_2-e_3-2e_4-e_5 & f_{34}f_{14}f_{25} & g_6f_{14}e_1 & e_6f_{45}g_5 & f_{24}f_{34}f_{15}\\ [1.8pt]
2\ell -e_1-e_2-e_4 & e_5f_{14}f_{25} & f_{14}f_{12}e_1 & e_6e_3 g_5 & f_{24}e_5f_{15}\\ [1.8pt]
2\ell -e_1-e_2-e_6 & e_5f_{16}f_{25} & f_{16}f_{12}e_1 & e_4e_3f_5 & f_{26}e_5f_{15}\\ [1.8pt]
4\ell -e_1-e_2-2e_3-2e_4-2e_5-e_6 &g_2f_{34}f_{25} & g_2g_6e_1 & f_{45}f_{35}g_5 & f_{34}g_1f_{15}\\ [1.8pt]
3\ell -e_1-e_3-e_4-e_5-2e_6 & f_{46}f_{13}f_{56} & e_2f_{56}g_5 & g_2f_{16}e_1 & f_{46}f_{36}f_{15}\\ [1.8pt]
4\ell -e_1-e_2-2e_3-e_4-2e_5-2e_6 & g_1f_{13}f_{56} & g_2g_4e_1 & f_{35}f_{56}g_5 & g_1f_{36}f_{15}\\ [1.8pt]
3\ell -e_1-2e_3-e_4-e_5-e_6 & f_{34}f_{13}f_{56} & e_2f_{35}g_5 & g_2f_{13}e_1 & f_{34}f_{36}f_{15}\\ [1.8pt]
2\ell -e_1-e_3-e_5 & e_5f_{13}f_{56} & e_2e_4g_5 & f_{16}f_{13}e_1 & e_5f_{36}f_{15}\\ [1.8pt]
3\ell -e_1-e_2-e_3-e_5-2e_6 & f_{26}f_{13}f_{56} & g_4f_{16}e_1 & e_4f_{56}g_5 & f_{26}f_{36}f_{15}\\ [1.8pt]
3\ell -e_1-e_3-2e_4-e_5-e_6 & f_{34}f_{14}f_{56} & e_2f_{45}g_5 & g_2f_{14}e_1 & f_{34}f_{46}f_{15}\\ [1.8pt]
2\ell -e_1-e_4-e_6 & e_5f_{14}f_{56} & e_2e_3g_5 & f_{16}f_{14}e_1 & e_5f_{46}f_{15}\\ [1.8pt]
2\ell -e_1-e_4-e_6 & f_{46}e_3f_{13} & e_2e_3g_5 & f_{16}f_{14}e_1 & e_5f_{46}f_{15}\\ [1.8pt]
2\ell -e_1-e_3-e_4 & f_{34}e_3f_{13} & e_2e_6g_5 & f_{14}f_{13}e_1 & e_5f_{34}f_{15}\\ [1.8pt]
2\ell -e_1-e_4-e_5 & f_{45}e_3f_{13} & e_6f_{14}f_{56} & f_{34}e_3f_{15} & f_{14}e_1f_{15}\\ [1.8pt]
2\ell -e_1-e_2-e_4 & f_{24}e_3f_{13} & f_{14}f_{12}e_1 & e_6e_3g_5 & f_{24}e_5f_{15}\\ [1.8pt]
2\ell -e_1-e_2-e_6 & f_{26}e_3f_{13} & f_{16}f_{12}e_1 & e_4e_3g_5 & f_{26}e_5f_{15}\\ [1.8pt]
3\ell -e_1-e_3-2e_4-e_5-e_6 & f_{45}f_{46}f_{13} & e_2f_{45}g_5 & g_2f_{14}e_1 & f_{34}f_{46}f_{15}\\ [1.8pt]
4\ell -e_1-e_2-2e_3-2e_4-2e_5-e_6 & f_{45}g_1f_{13} & g_2g_6e_1 & f_{45}f_{35}g_5 & f_{34}g_1f_{15}\\ [1.8pt]
3\ell -e_1-e_2-e_3-2e_4-e_5 & f_{24}f_{45}f_{13} & g_6f_{14}e_1 & e_6f_{45}g_5 & f_{24}f_{34}f_{15}\\ [1.8pt]
5\ell-2e_1-2e_2-2e_3-2e_4-2e_5-2e_6 & g_2f_{12}g_1 & g_6g_4f_{46} & g_2f_{25}g_5 & g_1g_5f_{15}\\ [1.8pt]
\hline 
\end{array}\]
and the cubic generator of degree $-K$ given by the first monomial in the table below. 
\[\begin{array}{|c|c|cccc|}
\hline
3\ell -e_1-e_2-e_3-e_4-e_5-e_6 & 
f_{46}f_{13}f_{25} & 
e_2f_{25}g_5 & f_{23}f_{46}f_{15} & g_2f_{12}e_1 & f_{26}f_{34}f_{15}\\
\hline
\end{array}\]
\begin{lemma} \label{lem:ugly1} For every Del Pezzo surface $X_6$, the monomials
\[\begin{array}{ccccc}
f_{34}f_{16}f_{25} & f_{46}f_{13}f_{25} & f_{23}f_{46}f_{15} & g_2f_{12}e_1 & f_{26}f_{34}f_{15}\\
\end{array}\]
are linearly dependent modulo $Q_6$.
\end{lemma}
\begin{proof} The ideal $Q_6$ contains elements 
\[\begin{array}{c}
p_1=f_{12}\bigl(a_1f_{46}f_{35}+a_2f_{34}f_{56}+a_3g_2e_1\bigr)\\[2.5pt]
p_2=f_{34}\bigl(b_1f_{12}f_{56}+b_2f_{16}f_{25}+b_3f_{26}f_{15}\bigr)\\[2.5pt]
p_3=f_{46}\bigl(c_1f_{12}f_{35}+c_2f_{13}f_{25}+c_3f_{23}f_{15}\bigr)\\[2.5pt]
\end{array}\]
with $a_i,b_i,c_i$ nonzero constants. Thus $p_1-\frac{a_1}{c_1}p_3$ is a nonzero polynomial of the form
\[d_1f_{34}f_{12}f_{56}+d_2f_{46}f_{13}f_{25}+d_3f_{23}f_{46}f_{15}+d_4g_2f_{12}e_1\]
with $d_1\neq 0$ since otherwise, evaluation at a point on the strict transform of the line through points 4 and 6 and in no other exceptional curve would show that all coefficients are zero. As a result $p_1-\frac{a_1}{c_1}p_3-\frac{b_1}{d_1}p_2\in Q_6$ is the required linear dependency relation.
\end{proof}
\begin{lemma}\label{lem:cubicgenerator} If $X_6$ is a Del Pezzo surface without Eckart points, the ideal $Q_6$ contains an element of the form
\[a_1f_{46}f_{13}f_{25}+a_2e_2f_{25}g_5+a_3f_{23}f_{46}f_{15}+a_4g_2f_{12}e_1+a_5f_{26}f_{34}f_{15}\]
where the $a_i$ are constants and $a_1\neq 0$.
\end{lemma}
\begin{proof} By Lemma~\ref{lem:ugly1} the ideal $Q_6$ contains a nonzero element of the form
\[s=b_1f_{34}f_{16}f_{25}+ b_2f_{46}f_{13}f_{25}+b_3f_{23}f_{46}f_{15}+b_4g_2f_{12}e_1+b_5f_{26}f_{34}f_{15}\]
We show that $b_1$ is nonzero if the surface has no Eckart points. Otherwise evaluation at the intersection points $f_{46}\cap f_{15}$ and $f_{46}\cap f_{12}$ shows that $b_4=b_5=0$ since no triple of exceptional curves has a common point. This forces $b_2=b_3=0$ yielding a contradiction.\\
Now  $Q_6$ also contains a relation of the form
\[r=f_{25}\bigl(c_1f_{34}f_{16}+c_2f_{46}f_{13}+c_3e_2g_5\bigr)\]
so $r-\frac{c_1}{b_1}s$ is a nonzero element of $Q_6$ of the desired form and the reasoning of the first paragraph (on the pairs of curves $f_{15},g_5$ and $f_{15},e_1$) shows that if $X_6$ has no Eckart points then $a_1\neq 0$. 
\end{proof}
The Table T4 contains the weight vectors $w$ which lead to all quadratic initial ideals of $C_5$ (up to the action of the Weyl Group $W_5=D_5$).
\[\begin{array}{c|c|c|c|c|c|c|c|c|c|c|c|c|c|c|c|}
\hline
T4\\
\hline
f_{12} & f_{13} & f_{23} &f_{14} & f_{24} & f_{34} & f_{15} & f_{25} & f_{35} & f_{45} & e_1 & e_2 & e_3 & e_4 & e_5 & g\\
25 & 21 &  11 &  19 &  13 & 9 &  13 & 9 & 1 &  11 &  10 &  10 &  10 &  10 &  10 &  10\\
18 &  15 &  12 &  14 &  10 & 9 &  12 &  10 &  11 &  5 &  10 &  10 &  10 &  10 &  10 &  10\\
25 &  19 &  15 &  17 &  11 & 9 &  15 &  13 &  17 &  3 &  10 &  10 &  10 &  10 &  10 &  10\\
17 &  15 &  13 &  14 &  11 &  10 &  12 &  11 &  12 &  5 &  10 &  10 &  10 &  10 &  10 &  10\\
19 &  12 &  15 &  11 &  13 &  10 &  10 &  14 & 8 & 5 &  10 &  10 &  10 &  10 &  10 &  10\\
20 &  12 &  15 & 9 &  11 & 4 &  11 &  16 &  10 &  10 &  10 &  10 &  10 &  10 &  10 &  10\\
39 & 31 &  23 &  25 &  21 & 19 &  27 &  21 & 17 & 9 &  20 &  20 &  20 &  20 &  20 &  20\\
17 &  15 &  10 &  14 &  11 & 9 &  12 &  10 & 6 &   12 &  10 &  10 &  10 &  10 &  10 &  10\\
27 &  19 &  17 &  17 &  11 & 9 &  17 & 13 &  19 &  3 &  10 &  10 &  10 &  10 &  10 &  10\\
18 &  17 & 13 &  15 &  12 &  11 &  11 &  10 & 8 & 11 &  10 &  10 &  10 &  10 &  10 &  10\\
18 &  16 &  13 &  14 &  12 &  11 &  11 &  10 & 7 & 11 &  10 &  10 &  10 &  10 &  10 &  10\\
20 &  14 &  16 &  10 &  11 & 6 &  12 &  15 &  11 &   11 &  10 &  10 &  10 &  10 &  10 &  10\\
18 &  16 &  12 &  14 &  11 &  10 &  11 & 9 & 6 &  10 &  10 &  10 &  10 &  10 &  10 &  10\\
25 & 21 &  13 &  19 &  15 &  11 &  13 &  11 & 3 &   13 &  10 &  10 &  10 &  10 &  10 &  10\\
18 &  15 &  13 &  14 &  11 &  10 &  12 &  11 &  12 &  6 &  10 &  10 &  10 &  10 &  10 &  10\\
20 &  13 &  16 & 9 &  11 & 5 &  11 &  15 &  10 &  10 &  10 &  10 &  10 &  10 &  10 &  10\\
20 &  13 &  15 &  10 &  11 & 5 &  12 &  16 &  11 &   11 &  10 &  10 &  10 &  10 &  10 &  10\\
21 & 27 &  15 &  17 &  11 &  13 &  13 &  13 &  11 &  1 &  10 &  10 &  10 &  10 &  10 &  10\\
\hline
\end{array}\]


\begin{thebibliography}{99}
%
\bibitem{COX}
Cox, D., {\it The Homogeneous Coordinate Ring of a Toric Variety}, 
J. Algebr. Geom. 4 {\bf 1}, 17-50, 1995.
%
\bibitem{GB}
Sturmfels, B., {\it Gr\"obner Basis and Convex Polytopes}, 
American Mathematical Society University Lecture Series, 1996.
%
\bibitem{DMc}
MacLagan, D., {\it Antichains of monomial ideals are finite}, 
Proc. Amer. Math. Soc. 129 {\bf 6}, 1609-1615, 2001.
%
\bibitem{HK}
Keel, S., Hu, Y., {\it Mori Dream Spaces}, 
Michigan Math Journal, 48, 331-348, 2000.
%
\bibitem{Man}
Manin, Yu., {\it Cubic Forms: Algebra, Geometry, Arithmetic}, 
North-Holland, 1974.
%
\bibitem{BP}
Batyev, V., Popov, O., {\it The Cox rings of a Del Pezzo surface}, 
Arithmetic of higher-dimensional algebraic varieties, B. Poonen 
and Y. Tschinkel (eds.), Progress in Math. Birkh\"auser, 226, 2004.
\end{thebibliography}
\end{document}